\documentclass[12pt]{amsart}
\usepackage{eurosym}
\usepackage{amssymb}
\usepackage{amsfonts}
\usepackage{subfigure}
\usepackage{graphicx,color}

\theoremstyle{plain}

\newtheorem{conjecture}{Conjecture}
\newtheorem{corollary}{Corollary}

\newtheorem{definition}{Definition}

\newtheorem{proposition}{Proposition}
\newtheorem{remark}{Remark}

\newtheorem{theorem}{Theorem}
\numberwithin{equation}{section}

\renewcommand{\S}{\mathcal{S}}

\begin{document}
\title[A large sieve priciple for modulation spaces]{Donoho-Logan large
sieve principles for modulation and polyanalytic Fock spaces}
\author{\medskip Lu\'{\i}s Daniel Abreu and Michael Speckbacher}
\address{Acoustics Research Institute, Austrian Academy of Sciences,\\
Wohllebengasse 12-14, 1040 Vienna, Austria}
\date{}

\begin{abstract}
We obtain estimates for the $L^{p}$-norm of the \emph{short-time Fourier
transform (STFT)} for functions in modulation spaces, providing information
about the concentration on a given subset of $\mathbb{R}^{2}$, leading to
deterministic guarantees for perfect reconstruction using convex
optimization methods. More precisely, we will obtain large sieve
inequalities of the Donoho-Logan type, but instead of localizing the signals
in regions $T\times W$ of the time-frequency plane using the Fourier
transform to intertwine time and frequency, we will localize the
representation of the signals in terms of the short-time Fourier transform
in sets $\Delta $ with arbitrary geometry. At the technical level, since
there is no proper analogue of Beurling's extremal function in the STFT
setting, we introduce a new method, which rests on a combination of an
argument similar to Schur's test with an extension of Seip's local
reproducing formula to general Hermite windows. When the windows are Hermite
functions, we obtain local reproducing formulas for polyanalytic Fock spaces
which lead to explicit large sieve constant estimates and, as a byproduct,
to a reconstruction formula for $f\in L^{2}(\mathbb{R})$ from its STFT
values on arbitrary discs. A discussion on optimality follows, along the
lines of Donoho-Stark paper on uncertainty principles and signal recovery.
We also consider the case of discrete Gabor systems, vector-valued STFT
transforms and rephrase the results in terms of the polyanalytic
Bargmann-Fock transforms.
\end{abstract}

\dedicatory{To the memory of Kurt G\"{o}del}
\maketitle


\section{Introduction}

The large sieve principle covers a number of far reaching analysis
techniques, mostly aimed at solving problems in analytic number theory, but
which have also found applications in a number of other mathematical fields,
like probability \cite{Kowalski}, numerical \cite{AuB} and signal analysis 
\cite[Theorem 7]{DonohoLogan,DonohoStark}, to name a few. The terminology
stems from its number theory origins, which can be traced back to the sieve
of Eratosthenes. In number theory, the large sieve principle is mostly
concerned with asymptotic averages of arithmetic functions on integers
constrained by congruences modulo sets of primes. A typical example of the
large sieve principle is the inequality for trigonometric polynomials:%
\begin{equation}
\sum_{l=1}^{R}\left\vert \sum_{k=m+1}^{m+n}a_{k}e^{2\pi ikx_{l}}\right\vert
^{2}\leq \Delta (n,\delta )\sum_{k=m+1}^{m+n}\left\vert a_{k}\right\vert
^{2},  \label{LS}
\end{equation}%
where the points $x_{1},...,x_{R}$ are $\delta $-separated mod $1$. By
choosing them as fractions $p/q$ with $\mbox{gcd}\left( p,q\right) =1$,
several applications in number theory follow \cite{Kowalski,Montgomery}.

According to inequality (\ref{LS}), only a small energy portion of the
trigonometric polynomial is concentrated at the points $\alpha
_{1},...,\alpha _{R}$, with the constant $\Delta (n,\delta )$ controlling
the size of the fraction. Since the energy concentration inside a domain is
important to find optimal approximation methods required in signal recovery,
such an observation suggests applications in signal analysis, where one can
find a rich setting. Donoho and Logan \cite{DonohoLogan}, in a paper that,
together with \cite{DonohoStark}, spearheaded the modern theory of
Compressed Sensing \cite{CRT,Donoho,CS}, introduced the concept of \emph{%
maximum Nyquist density}, $\rho (T,W)$, which measures the sparsity of a
real band-limited signal on the time domain $T\subset \mathbb{R}$ with
band-size $W$: 
\begin{equation}
\rho (T,W):=W\cdot \sup_{t\in \mathbb{R}}\left\vert T\cap \lbrack
t,t+1/W]\right\vert \leq W\cdot \left\vert T\right\vert \text{.}  \label{MND}
\end{equation}%
Let us motivate our work with some results from \cite%
{DonohoLogan,DonohoStark}. If the set $T$ {has large area but small Lebesgue
measure on any interval of length $1/W$,} then $\rho (T,W)$\ can be
considerably small compared to the natural Nyquist density $W\cdot
\left\vert T\right\vert $. We will call such sets $T$ \emph{sparse in the
sense of Lebesgue measure}. Throughout the paper we will write $P_{T}f:=\chi
_{T}f$ to denote multiplication by the indicator function of $T$.

While analytic number theory is mostly interested in Hilbert space large
sieve inequalities, in signal analysis one finds remarkable applications of
Banach space large sieve inequalities, with a special emphasis on $L^{1}$%
-normed spaces. In \cite[Theorem 7]{DonohoLogan}, the authors considered the
space 
\begin{equation}
B_{1}(W):=\left\{ f\in L^{1}(\mathbb{R}):\ supp(\hat{f})\subseteq \lbrack
-\pi W,\pi W]\right\}  \label{B1}
\end{equation}%
and proved, for $\theta <\frac{2\pi }{W}$, the inequality%
\begin{equation}
\Vert P_{T}f\Vert _{1}\leq \frac{\pi W/2}{\sin (\pi W\theta /2)}\left(
\sup_{t\in \mathbb{R}}\left\vert T\cap \lbrack t,t+\theta ]\right\vert
\right) \cdot \Vert f\Vert _{1}\text{.}  \label{DLl1}
\end{equation}%
An inequality similar as (\ref{DLl1}) also holds if $\chi _{T}(x)dx$ is
replaced by {a positive $\sigma $-finite measure $\mu $}. It provides the
concentration bound%
\begin{equation*}
\delta _{1}(T):=\sup_{f\in B_{1}(W)}\frac{\Vert P_{T}f\Vert _{1}}{\Vert
f\Vert _{1}}\leq \frac{\pi }{2}\rho (T,W)\text{,}
\end{equation*}%
resulting in sufficient conditions for perfect reconstruction of a
band-limited signal corrupted by sparse noise using $L^{1}$-norm
minimization.

In this paper, we will obtain large sieve inequalities of the Donoho-Logan
type, but instead of localizing the signals in regions $T\times W$ of the
time-frequency plane using the Fourier transform to intertwine time and
frequency, we will localize the representation of the signals in terms of the%
{\ \emph{short-time Fourier transform (STFT)}} {\ 
\begin{equation}
V_{g}f(x,\xi )=\int_{%
%TCIMACRO{\U{211d} }%
%BeginExpansion
\mathbb{R}
%EndExpansion
}f(t)\overline{g(t-x)}e^{-2\pi i\xi t}dt\text{,}  \label{Gabor}
\end{equation}%
}to general regions $\Delta $ of the time-frequency plane. Instead of {t}he
maximum Nyquist density (\ref{MND}) we will use the following concept of
planar maximum Nyquist density $\rho (\Delta ,R)$ introduced in \cite%
{AbreuSpeck}: 
\begin{equation}
\rho (\Delta ,R):=\sup_{z\in \mathbb{R}^{2}}|\Delta \cap (z+D_{1/R})|\leq
|\Delta |\text{,}  \label{planarNyquist}
\end{equation}%
where $D_{1/R}\subset \mathbb{R}^{2}$ is the disc of radius $1/R$ centered
in the origin. If the set $\Delta $ is sparse in the sense of Lebesgue
measure (small concentration in any disc of radius $1/R$), then $\rho
(\Delta ,R)$\ can be considerably smaller than the natural Nyquist density\ $%
|\Delta |$ (see \cite{Daubechies,defeno02,Sampling,Ly,Sampling1,Seip0} for
natural Nyquist densities in the context of Fock and Gabor spaces). Instead
of measuring the concentration of band-limited signals in a time-limited
region $T\subset \mathbb{R}$, we will measure the joint time-frequency
content on a region $\Delta \subset \mathbb{R}^{2}$. In its most general
version, our results can be seen as estimates on the bounds of Bessel
measures for the STFT \cite{Ascensi,OrtegaCerda}. By selecting Hermite
functions as windows in (\ref{Gabor}), good explicit estimates in terms of $%
\rho (\Delta ,R)$\ can be obtained.\ The following is a sample of our
findings in the $L^{1}$-case (we will prove it for $1\leq p<\infty $). The
modulation space $M^{1}$, also known as Feichtinger's algebra $S_{0}$\ \cite%
{Fei},\ will play the role of the space $B_{1}(W)$\ in \cite{DonohoLogan}.

\noindent \textbf{Theorem. }\emph{Let }$\Delta \subset \mathbb{R}^{2}$\emph{%
\ be measurable and }$f\in M^{1}$\emph{. Denote by }$h_{r}$\emph{\ the }$rth$
\emph{Hermite function}. \emph{For every }$0<R<\infty $\emph{,}%
\begin{equation}
\Vert V_{h_{r}}f\cdot \chi _{\Delta }\Vert _{1}\leq \frac{\rho (\Delta ,R)}{%
C_{r}(R)}\Vert V_{h_{r}}f\Vert _{1}\text{,}\ 
\end{equation}%
\emph{where the constant }$C_{r}(R)$\emph{\ is explicitly determined.}

This will provide estimates for the concentration bound%
\begin{equation*}
\delta (\Delta ):=\sup_{f\in M^{1}}\dfrac{\Vert V_{h_{r}}f\cdot \chi
_{\Delta }\Vert _{1}}{\Vert V_{h_{r}}f\Vert _{1}}\text{,}
\end{equation*}%
and, consequently, conditions for perfect reconstruction of a $M^{1}$
function corrupted by sparse noise using $L^{1}$-norm minimization (precise
statements are given in section 5.1).

The techniques of proof are new in sieve theory and, in particular,
different from those in \cite{DonohoLogan,Montgomery}, where Beurling's
extremal function \cite{Beurling} plays a key role. Since we are not aware
of a proper analogue of extremal function theory in the STFT setting, we had
to develop new methods, which essentially depend on combining an argument
similar to Schur's test with an extension of Seip's local reproducing
formula's \cite{Seip0} to general Hermite windows.

The paper is organized as follows. In the preliminaries section we gather
the essential background on time-frequency analysis. In Section \ref%
{sec:general-LS}, we formulate our results in a general Banach space
setting, provide a general discussion about the signal recovering
applications that motivate the results and extend Selberg-Bombieri
inequality \cite{Bombieri} to the continuous setting. In Section \ref%
{sec:local-repr} we restrict to Hermite windows. Then, we extend Seip's
local reproducing formulas \cite{Seip0} to polyanalytic Fock spaces
associated with the Landau levels and use them to obtain explicit estimates
for the maximum Nyquist density in the Hermite window case. Section \ref%
{sec:optimality} contains a discussion on optimality of the constants,
including a phase-space versions of a theorem by Donoho-Stark \cite[Theorem
10]{DonohoStark}. Moreover, we revisit the signal recovery problem in the
context of reconstructing STFT data. In Section 6 we obtain large sieve
inequalities for discrete and vector valued Gabor systems. Finally, we
conclude with a section discussing some open problems.

\section{Preliminaries}

\subsection{The short-time Fourier transform}

Let $z=(x,\xi ),\ w=(y,\eta )\in \mathbb{R}^{2}$ and $g\in L^{2}(\mathbb{R})$%
. A time-frequency shift of the function $g$ is defined as 
\begin{equation}
\pi (z)g(t):=M_{\xi }T_{x}g(t)=e^{2\pi it\xi }g(t-x),
\end{equation}%
where $T_{x}$ denotes the translation operator and $M_{\xi }$ the modulation
operator. The composition of two time-frequency shifts is given by 
\begin{equation}
\pi (z)\pi (w)=e^{-2\pi ix\eta }\pi (z+w)  \label{compos-tf}
\end{equation}%
and the adjoint operator of $\pi (z)$ is 
\begin{equation}
\pi (z)^{\ast }=e^{-2\pi ix\xi }\pi (-z).  \label{adjoint-tf}
\end{equation}%
The short-time Fourier transform (STFT) or Gabor transform of a function $f$
with window $g$ is defined by 
\begin{equation*}
V_{g}f(z):=\langle f,\pi (z)g\rangle =\int_{%
%TCIMACRO{\U{211d} }%
%BeginExpansion
\mathbb{R}
%EndExpansion
}f(t)\overline{g(t-x)}e^{-2\pi i\xi t}dt.
\end{equation*}%
An important property of the STFT is the so called orthogonality relation 
\begin{equation}
\langle V_{g_{1}}f_{1},V_{g_{2}}f_{2}\rangle =\langle f_{1},f_{2}\rangle
\langle g_{2},g_{1}\rangle .  \label{ortho-rel}
\end{equation}%
In particular, if $\Vert g\Vert _{2}=1$, then 
\begin{equation*}
\Vert V_{g}f\Vert _{2}=\Vert f\Vert _{2}
\end{equation*}%
and $V_{g}:L^{2}(\mathbb{R})\rightarrow L^{2}(\mathbb{R}^{2})$ is an
isometry mapping onto the reproducing kernel Hilbert space 
\begin{equation*}
\mathcal{V}_{g}:=\{V_{g}f:\ f\in L^{2}(\mathbb{R})\}\subset L^{2}(\mathbb{R}%
^{2}).
\end{equation*}%
The corresponding reproducing equation is 
\begin{equation*}
V_{g}f(z)=\langle V_{g}f,K_{g}(z,\cdot )\rangle \text{,}
\end{equation*}%
where $K_{g}(z,w):=\langle \pi (w)g,\pi (z)g\rangle $. Thus, one can write
the orthogonal projection of $L^{2}(\mathbb{R}^{2})$ on $\mathcal{V}_{g}$ as
the integral operator with kernel $K_{g}$. We note in passing that $%
K_{g}(z,w)$ is the correlation kernel in the Weyl-Heisenberg ensemble \cite%
{WH,abgrro17}.

%Consider then a function $f:{\mathbb{R}^{1}}\rightarrow {\mathbb{C}}$, and let \begin{equation}{H_{\Delta }}f(t)=\int_{{\mathbb{R}^{2}}}1_{\Delta }(x,\xi )V_{g}f(x,\xi)g(t-x)e^{2\pi i\xi t}dxd\xi \text{,}\qquad t\in {\mathbb{R}}\text{.}\label{localization}\end{equation}%The indicator function $1_{\Delta }(x,\xi )$ is called the symbol of ${H_{\Delta }}$.

\subsection{Hermite functions and complex Hermite polynomials}

In time-frequency analysis, a particular interest is given to functions well
concentrated in both time and frequency. A class of such functions is given
by the Hermite functions $h_{r}$ defined as 
\begin{equation*}
h_{r}(t)=\frac{2^{1/4}}{\sqrt{r!}}\left( \frac{-1}{2\sqrt{\pi }}\right)
^{r}e^{\pi t^{2}}\frac{d^{r}}{dt^{r}}\left( e^{-2\pi t^{2}}\right) ,\qquad
r\geq 0\text{.}
\end{equation*}%
The collection $\{h_{r}\}_{r\geq 0}$ forms an orthonormal basis for $L^{2}(%
\mathbb{R})$, minimizes the uncertainty principle \cite{Jaming} and
optimizes the joint time-frequency concentration on discs \cite%
{Daubechies,Seip0}. In \cite{CharlyYurasuper}, precise lattice conditions
for vector valued frames with Hermite function have been obtained, which,
combined with Vasilevski%
%TCIMACRO{\U{b4}}%
%BeginExpansion
\'{}%
%EndExpansion
s work \cite{VasiFock}, inspired the study of sampling and interpolation
problems in Fock spaces of polyanalytic functions \cite{Abreusampling}, a
hierarchy of function spaces which seems to be ubiquitous in several
mathematical models \cite{abfei14}. The Hermite function $h_{0}$ is the
Gaussian function explicitly given by 
\begin{equation*}
\varphi (t)=h_{0}(t)=2^{1/4}e^{-\pi t^{2}}\text{.}
\end{equation*}%
We will also use the so-called complex Hermite polynomials \cite%
{Ghanmi,Ismail}: 
%\begin{equation}H_{j,r}(z,\overline{z})=\left\{ \begin{tabular}{l}${\sqrt{\frac{r!}{j!}}\pi ^{\frac{j-r}{2}}z^{j-r}L_{r}^{j-r}\left( \pi\left\vert z\right\vert ^{2}\right) ,\qquad j>r\geq 0,}$ \\ ${\left( -1\right) ^{r-j}\sqrt{\frac{j!}{r!}}\pi ^{\frac{r-j}{2}}\overline{z}^{r-j}L_{j}^{r-j}\left( \pi \left\vert z\right\vert ^{2}\right) ,\qquad 0\leq j\leq r}$,\end{tabular}\ \right.   \label{ComplexHermite}\end{equation}
\begin{equation}
H_{j,r}(z,\overline{z})=\left\{ 
\begin{tabular}{l}
${\sqrt{\frac{r!}{j!}}\pi ^{\frac{j-r}{2}}z^{j-r}L_{r}^{j-r}\left( \pi
\left\vert z\right\vert ^{2}\right) ,\qquad j>r\geq 0}${,} \\ 
${\left( -1\right) ^{r-j}\sqrt{\frac{j!}{r!}}\pi ^{\frac{r-j}{2}}\overline{z}%
^{r-j}L_{j}^{r-j}\left( \pi \left\vert z\right\vert ^{2}\right) ,\qquad
0\leq j\leq r}$,%
\end{tabular}%
\ \right.  \label{ComplexHermite}
\end{equation}%
where ${L_{r}^{j-r}}$ stands for the generalized Laguerre polynomials
defined via the recurrence relation $L_{0}^{\alpha }(x)=1$, $L_{1}^{\alpha
}(x)=1+\alpha -x$ and 
\begin{equation*}
L_{j+1}^{\alpha }(x)=\frac{2j+1+\alpha -x}{j+1}L_{j}^{\alpha }(x)-\frac{%
j+\alpha }{j+1}L_{j-1}^{\alpha }(x),\qquad j\geq 1.
\end{equation*}%
If $\alpha \geq 0$, then $L_{j}^{\alpha }$ has the following closed form 
\begin{equation*}
L_{j}^{\alpha }(x)=\sum\limits_{i=0}^{j}(-1)^{i}\binom{j+\alpha }{j-i}\frac{%
x^{i}}{i!},\qquad x\in {\mathbb{R}},\qquad j,\alpha \in \mathbb{N}_{0}.
\end{equation*}%
Complex Hermite polynomials satisfy the doubly-indexed orthogonality 
\begin{equation*}
\int_{\mathbb{C}}H_{j,r}(z,\overline{z})\overline{H_{j\prime ,r\prime }(z,%
\overline{z})}e^{-\pi \left\vert z\right\vert ^{2}}dz={\delta }_{jj\prime }{%
\delta }_{rr\prime },
\end{equation*}%
and provide a basis for the space $L^{2}\big(\mathbb{C},e^{-\pi \left\vert
z\right\vert ^{2}}\big)$ \cite{abgroe12,HendHaimi}. The relation between
time-frequency analysis and polyanalytic functions \cite%
{Abreustructure,Abreusampling} can be understood in terms of the \emph{%
Laguerre connection} \cite[Chapter 1.9]{Folland} 
\begin{equation}
V_{h_{r}}h_{j}(x,-\xi )=e^{i\pi x\xi -\tfrac{\pi }{2}\left\vert z\right\vert
^{2}}H_{j,r}(z,\bar{z})\text{.}  \label{l9b}
\end{equation}%
The closed form of the reproducing kernels $K_{h_{r}}$ reads \cite%
{Abreustructure} 
\begin{equation*}
K_{h_{r}}(z,w)=\langle \pi (w)h_{r},\pi (z)h_{r}\rangle =e^{i\pi (x+y)(\xi
-\eta )}L_{r}^{0}(\pi |z-w|^{2})e^{-\pi |z-w|^{2}/2}\text{.}
\end{equation*}%
Consequently, 
\begin{equation}
|K_{h_{r}}(z,w)|=L_{r}^{0}(\pi |z-w|^{2})e^{-\pi |z-w|^{2}/2}.
\label{abs-kernel}
\end{equation}%
The kernel $K_{h_{r}}$ describes the orthogonal projection onto the
Bargmann-Fock space of pure polyanalytic functions of type $r$ (see Remark %
\ref{remark-poly}), which is precisely the $rth$-eigenspace of the Euclidean
Landau operator with a constant magnetic field \cite%
{AoP,Zouhair,HendHaimi,Ruz}.

\subsection{Modulation spaces and (Banach-)frame theory}

In order to quantitatively measure the behavior of a class of
transformations generated by an integrable group representation, Feichtinger
and Gr{\"{o}}chenig developed the theory of coorbit spaces \cite%
{fegr88,fegr89,fegr89II}. We will consider the particular instance of the
coorbit spaces associated with the short-time Fourier transform, the so
called modulation spaces \cite{Fei-Mod}. Let $g$ be a window function
satisfying $\Vert V_{g}g\Vert _{1}<\infty $. The modulation space $M^{p}$ is
defined as 
\begin{equation}
M^{p}:=\{f\in \S ^{\prime }(\mathbb{R}):\ V_{g}f\in L^{p}(\mathbb{R}^{2})\},
\label{def-modul-space}
\end{equation}%
where $\mathcal{S}^{\prime }(\mathbb{R})$ denotes the space of tempered
distributions. The space $M^{1}$, also known as Feichtinger's algebra $S_{0}$%
,\ will play the role of $B_{1}(W)$\ in \cite{DonohoLogan} as the
fundamental space for applications in signal recovery using $L_{1}$%
-minimization. Since the reproducing kernel property extends to $M^{p}$, the
range of the short-time Fourier transform on $M^{p}$ can be characterized in
terms of projections on $L^{p}(\mathbb{R}^{2})$ as follows:%
\begin{equation}
\mathcal{V}_{g}^{p}:=V_{g}(M^{p})=\{F\in L^{p}(\mathbb{R}^{2}):\ \langle
F,K_{g}(z,\cdot )\rangle =F(z)\}.  \label{corr-princ}
\end{equation}%
Let $\Lambda \subset \mathbb{R}^{2}$ be a discrete set. The family $\{\pi
(\lambda )g\}_{\lambda \in \Lambda }$ is a Gabor frame for $L^{2}(\mathbb{R}%
) $ if there exist positive constants $A,B$ such that 
\begin{equation}
A\Vert f\Vert _{2}^{2}\leq \sum_{\lambda \in \Lambda }|\langle f,\pi
(\lambda )g\rangle |^{2}\leq B\Vert f\Vert _{2}^{2},\qquad \forall f\in
L^{2}(\mathbb{R}).  \label{gabor-frame}
\end{equation}%
It follows from the theory of modulation spaces that a Gabor frame $%
G(g,\Lambda )$ with $g\in M^{1}$ is also a Gabor frame for every space $%
M^{p} $ with $p\geq 1$ \cite[Theorem 13.6.1]{Charly}. This is summarized in
the following remark.

\begin{remark}
\label{banach-gabor-frame} Assume that $g\in M^{1}$ and $G(g,\Lambda )$ is a
Gabor frame for $L^{2}(\mathbb{R})$. There exist two constants $A^{\prime
},B^{\prime }>0$ such that 
\begin{equation*}
A^{\prime }\Vert f\Vert _{M^{p}}^{p}\leq \sum_{\lambda \in \Lambda }|\langle
f,\pi (\lambda )g\rangle |^{p}\leq B^{\prime }\Vert f\Vert
_{M^{p}}^{p},\qquad \forall f\in M^{p}.
\end{equation*}%
Furthermore, there exists a dual window $\gamma \in M^{1}$ such that 
\begin{equation*}
f=\sum_{\lambda \in \Lambda }\langle f,\pi (\lambda )g\rangle \pi (\lambda
)\gamma ,\qquad \forall f\in M^{p},
\end{equation*}%
with unconditional convergence in $M^{p}$ if $p<\infty $ and weak-$\ast $
convergence in $M^{\infty }$.
\end{remark}

\section{A general large sieve principle}

\label{sec:general-LS}

In this section we will first formulate our problem and applications in a
very general setting. Later, particular situations with more structure will
be considered, where explicit estimates for the large sieve constant can be
obtained.

\subsection{Sieving inequalities in Banach spaces}

Let $(X,\mu )$, $(X,\nu)$ be measure spaces and $B_{p}\subset L^{p}(X,\nu)$
a Banach space. We derive a bound on the embedding $(B_{p},\Vert \cdot \Vert
_{L^p_\nu})\hookrightarrow (B_{p},\Vert \cdot \Vert _{L_{\mu }^{p}})$ using
an argument similar to Schur's test. Define the Banach algebra $\mathcal{A}%
_{\mu }$ of hermitian integral kernels, $K(x,y)=\overline{K(y,x)}$, equiped
with the norm 
\begin{equation*}
\Vert K\Vert _{\mathcal{A}_{\mu }}:=\sup_{y\in X}\int_{X}|K(x,y)|d\mu (x)
\end{equation*}%
and multiplication rule 
\begin{equation*}
K_{1}\circ K_{2}(x,y)=\int_{X}K_{1}(x,z)K_{2}(z,y)d\mu (z).
\end{equation*}

\begin{proposition}
\label{prop-measure-bound} Let $\mu $ be a positive $\sigma $-finite measure
on $X$, $B_{1}\subset L^{1}(X,\nu)$ and $K\in \mathcal{A}_{\mu }$ be 
%$\varepsilon$-concentrated w.r.t. $\Omega_R$  
such that $\mathcal{K}:B_{1}\rightarrow B_{1}$, 
\begin{equation*}
\mathcal{K}F(x):=\int_{X}F(y)K(x,y)d\nu(y),
\end{equation*}%
is bounded and boundedly invertible on $B_{1}$. Then, for every $F\in B_{1}$%
, defined as 
\begin{equation}
\frac{\int_{X}|F|d\mu }{\int_X |F|d\nu}\leq \theta (K)\cdot \Vert K\Vert _{%
\mathcal{A}_{\mu }},  \label{eq-norm-bound2}
\end{equation}%
where 
\begin{equation*}
\theta (K):=\sup_{\phi \in B_{1}}\left( \dfrac{\Vert \phi \Vert _{L^1_\nu}}{%
\Vert \mathcal{K}\phi \Vert _{L^1_\nu}}\right) .
\end{equation*}%
%
%
%
%
%
%
%
%
%
%
%
%
%
%
%
%
%
%
%
%
%
%
%
%
%
% and
% $$
%\Lambda_m(\mu,K,R):=\sup_{x\in\R^d}\left(\int_{\R^d} \chi_{D_R}( x-y)|K(x,y)|\frac{m(y)}{m(x)}d\mu(y)\right).
%$$
\end{proposition}

\textbf{Proof:\ } Let $F^{\ast }\in B_{1}$ be the unique function such that $%
\mathcal{K}F^{\ast }=F$. Using Fubini's theorem, we have 
\begin{align*}
\int_{X}|F(x)|d\mu (x) &\leq \int_{X}\int_{X}|F^{\ast }(y)K(x,y)|d\nu(y)d\mu
(x) \\
&=\int_{X}|F^{\ast }(y)|\int_{X}|K(x,y)|d\mu (x)d\nu(y) \\
&\leq \Vert F^{\ast }\Vert _{L^1_\nu}\cdot \Vert K\Vert _{\mathcal{A}_{\mu }}
\\
&=\dfrac{\Vert F^{\ast }\Vert _{L^1_\nu}}{\Vert \mathcal{K}F^{\ast }\Vert
_{1}}\cdot \Vert K\Vert _{\mathcal{A}_{\mu }}\cdot \Vert F\Vert _{L^1_\nu} \\
&\leq \theta (K)\cdot \Vert K\Vert _{\mathcal{A}_{\mu }}\cdot \Vert F\Vert
_{L^1_\nu}\text{.}
\end{align*}%
\hfill $\Box $\newline
The following result follows immediately by complex interpolation.

\begin{corollary}
\label{measure-coroll} Let $\Vert F\Vert _{L_{\mu }^{\infty }}\leq C_{\infty
}\Vert F\Vert _{L^\infty_\nu },\ \forall F\in B_{\infty }$ and $1\leq
p<\infty $. Under the assumptions of Proposition \ref{prop-measure-bound},
the following inequality holds 
\begin{equation}
\frac{\int_{X}|F|^{p}d\mu }{\int_X|F|^{p}d\nu}\leq C_{\infty }^{p-1}\cdot
\theta (K)\cdot \Vert K\Vert _{\mathcal{A}_{\mu }},\qquad \forall F\in
B_{1}\cap B_{\infty }.
\end{equation}
\end{corollary}

\subsection{Localization and signal recovery}

We now explain how the general estimates of the previous section can provide
useful information for two signal recovery scenarios. Let $B_{p}\subset
L^{p}(X)$ be a Banach space, where $\Delta\subset X\subset \mathbb{R}^{d}$
and define 
\begin{equation*}
\delta_p (\Delta ):=\sup_{F\in B_{p}}\dfrac{\int_{\Delta }|F(x)|^pdx}{%
\int_{X}|F(x)|^pdx}.
\end{equation*}

\subsubsection{Scenario 1: Perfect recovery of a signal corrupted by sparse
noise by $L^{1}$-minimization}

Let us assume that we observe a noisy version of a signal $F\in B_1$. In
addition let us assume that the noise $N$ has arbitrary but finite $L^1$%
-norm and is supported on an unknown set $\Delta$. The following recovery
result of Donoho and Stark \cite{DonohoStark} tells us that perfect
reconstruction is possible if only less than $1/2$ of the mass of a function
in $B_1$ can be concentrated on $\Delta$.

\begin{proposition}
\label{l1-sparse} Let $\delta_1(\Delta)<1/2$ and $G=F+N$, where $F\in B_1$
and $supp(N)\subset \Delta$. Then perfect reconstruction of $F$ is possible
via 
\begin{equation*}
F=\arg\min_{B\in B_1} \|B-G\|_1.
\end{equation*}
\end{proposition}

\subsubsection{Scenario 2: Approximated recovery of missing data by $L^{1}$%
-minimization}

Let us now assume that the information of the signal on $\Delta$ is missing,
so that one observes 
\begin{equation*}
H(x)=\left\{%
\begin{array}{ll}
(F+N)(x), & \mbox{for }x\notin \Delta \\ 
0, & \mbox{for }x\in \Delta%
\end{array}%
\right. ,
\end{equation*}
with the noise $N$ having small norm. The task of finding approximations of $%
F$ from $G$ is also known as the inpainting problem in signal processing 
\cite{inpaint1,inpaint2}.

\begin{proposition}
\label{l1-missing} Let $F\in B_1$, $\|N\|_1\leq\varepsilon$ and $%
\delta_1(\Delta)<1$. Set 
\begin{equation*}
\beta(H):=\arg\min_{S\in B_1}\|(I-P_\Delta)(H-S)\|_1.
\end{equation*}
For any solution $\beta(H)$ we have 
\begin{equation*}
\|F-\beta(H)\|_1\leq \frac{2\varepsilon}{1-\delta_1(\Delta)}.
\end{equation*}
\end{proposition}

\textbf{Proof:\ } The proof follows from the estimate 
\begin{equation*}
\Vert (I-P_{\Delta })(H-\beta (H))\Vert _{1}\leq \Vert (I-P_{\Delta
})(F+N-F)\Vert _{1}\leq \Vert n\Vert _{1}\leq \varepsilon \text{,}
\end{equation*}%
since 
\begin{align*}
\Vert F-\beta (H)\Vert _{1} &=\Vert (I-P_{\Delta })(F-\beta (H))\Vert
_{1}+\Vert P_{\Delta }(F-\beta (H))\Vert _{1} \\
&\leq \Vert (I-P_{\Delta })(H-\beta (H))\Vert _{1}+\Vert (I-P_{\Delta
})(F-H)\Vert _{1}+\delta (\Delta )\Vert F-\beta (H)\Vert _{1} \\
&\leq 2\varepsilon +\delta _{1}(\Delta )\Vert F-\beta (H)\Vert _{1}\text{.}
\end{align*}%
\hfill $\Box $

\begin{remark}
It is essential to assume the existence of a solution to the minimization
problem of Proposition \ref{l1-missing}. Moreover, the solution is not
necessarily unique.
\end{remark}

A similar result can be shown using $L^{2}$-minimization:

\begin{proposition}
Let $F\in B_{2}$ and define 
\begin{equation*}
\gamma (H):=\text{arg}\min_{S\in B_{2}}\Vert (I-P_{\Delta })(H-S)\Vert _{2}.
\end{equation*}%
If $\delta _{2}(\Delta )<1,$ then for any solution $\gamma (H)$ it holds 
\begin{equation*}
\Vert F-\gamma (H)\Vert _{2}^{2}\leq \frac{4\varepsilon ^{2}}{1-\delta
_{2}(\Delta )}.
\end{equation*}
\end{proposition}

\textbf{Proof:\ } Use the argument of Proposition \ref{l1-missing} and the
inequality $(a+b)^{2}\leq 2a^{2}+2b^{2}$.\hfill $\Box $

\subsection{Selberg's inequality}

A similar estimate for the case $p=2$ can be formulated as a direct
corollary of Selberg's inequality, as stated by Bombieri in \cite{Bombieri}.
Since the proof follows the same structure of the proof in \cite[Proposition
1]{Bombieri} by replacing sums by integrals, we leave the details to the
interested reader.

\begin{proposition}
Let $\mathcal{H}$ be a separable Hilbert space, $(X,\mu )$ a measure space, $%
\Delta \subseteq X$ and ${\psi :X\rightarrow }\mathcal{H}$. Then 
\begin{equation*}
\int_{\Delta }\frac{\left\vert \left\langle f,{\psi }_{x}\right\rangle
\right\vert ^{2}}{\int_{\Delta }\left\vert \left\langle {\psi }_{x},{\psi }%
_{y}\right\rangle \right\vert d\mu(y)}d\mu(x)\leq \Vert f\Vert ^{2},\qquad
\forall f\in \mathcal{H}\text{.}
\end{equation*}
\end{proposition}

\begin{corollary}
Under the same assumption as before, 
\begin{equation}
\frac{\int_{\Delta }|\langle {f},{\psi _{x}}\rangle |^{2}d\mu(x)}{\Vert
f\Vert ^{2}}\leq \sup_{x\in \Delta }\int_{\Delta }|\langle {\psi _{x}},{\psi
_{y}}\rangle |d\mu(y),\qquad \forall f\in \mathcal{H}\text{.}
\label{selberg1}
\end{equation}
\end{corollary}

If we take the $\psi_x$ to be time-frequency shifts of a function $g\in M^1$
with $\|g\|_2=1$, we obtain the following result.

\begin{corollary}
Let $\Delta \subset \mathbb{R}^{2}$ and $g\in M^{1}$. Then, for every $f\in
L^{2}(\mathbb{R})$, 
\begin{equation}
\frac{\int_{\Delta }|V_{g}f(z)|^{2}dz}{\Vert V_{g}f\Vert _{2}^{2}}\leq
\sup_{z\in \Delta }\int_{\Delta }|K_{g}(z,w)|dw\text{.}  \label{L2-stft}
\end{equation}
\end{corollary}

Note that this result yields a slight improvement to Proposition \ref%
{loc-general-g} of the next section in the case $p=2$, since the supremum
only needs to be taken over $\Delta $ and not over $\mathbb{R}^{2}$.

In the next sections, we will restrict the windows $g$ to the family of
Hermite functions, where, besides a sharp off diagonal fall of the kernel,
one can perform explicit computations and obtain workable explicit large
sieve constants.

\section{Local reproducing formulas and explicit maximum Nyquist density
estimates}

\label{sec:local-repr}

\subsection{Sieving inequalities for general short-time Fourier transforms}

Specializing the result for general Banach spaces from Corollary \ref%
{measure-coroll}, the following inequality for the STFT in modulation spaces
follows.

\begin{proposition}
\label{loc-general-g} Let $g\in M^{1}$ with $\Vert g\Vert _{2}=1$. For $f\in
M^{p}$, $1\leq p<\infty $, it holds 
\begin{equation}
\dfrac{\Vert V_{g}f\cdot \chi _{\Delta }\Vert _{p}^{p}}{\Vert V_{g}f\Vert
_{p}^{p}}\leq \sup_{z\in \mathbb{R}^{2}}\int_{\Delta }|K_{g}(z,w)|dw\text{.}
\label{prop-general-g-eq}
\end{equation}
\end{proposition}

Despite the apparent simplicity of the above inequality, it is virtually of
no use without information about the kernel $K_{g}(z,w)$. If the kernel has
proper off-diagonal decay properties, one expects to obtain good large sieve
constants and \eqref{prop-general-g-eq} can be simplified. Indeed, assume
that $K_{g}$ is $\varepsilon $-concentrated on $\Omega \subset \mathbb{R}%
^{2} $, more precisely, that 
\begin{equation*}
\sup_{z\in \mathbb{R}^{2}}\int_{\mathbb{R}^{2}\backslash (z+\Omega
)}|K_{g}(z,w)|dw<\varepsilon
\end{equation*}%
and $\Vert g\Vert _{2}=1$. Then, 
\begin{equation}
\dfrac{\Vert V_{g}f\cdot \chi _{\Delta }\Vert _{p}^{p}}{\Vert V_{g}f\Vert
_{p}^{p}}\leq \sup_{z\in \mathbb{R}^{2}}|\Delta \cap (z+\Omega
)|+\varepsilon \text{.}  \label{prop-general-g-eq2}
\end{equation}%
Depending on the window $g$, the set $\Omega $ may have to be chosen to be
big, leading to bad estimates in \eqref{prop-general-g-eq2}. Using the local
reproducing formulas of the next subsection, we will see that, for circular
domains $\Omega \subset \mathbb{R}^{2}$ and choosing windows from the
Hermite function sequence, one can obtain telling explicit estimates.

\subsection{Local reproducing formulas for Hermite windows and circular
domains}

With a slight abuse of language, double orthogonality often refers to
orthogonality in concentric domains. This is known to be the case for STFT's
of Hermite functions with Gaussian windows \cite{Seip0} which span the
Bargmann-Fock space of entire functions. In this contribution we show that
also the Hermite functions allow for local reproducing formulas, extending
the results in \cite[Proposition 4.2]{abgrro17} to Bargmann-Fock spaces of
polyanalytic functions. In fact we show even more: the reproducing kernel
corresponding to the Hermite function $h_{r}$ locally reproduces the
short-time Fourier transform using the window function $h_{j}$. At first
sight, this may be perceived as a counter-intuitive result, since 
\begin{equation*}
\int_{\mathbb{R}^{2}}V_{h_{r}}f(z)\overline{V_{h_{j}}h_{j}(z)}dz=0,\ %
\mbox{for }r\neq j,
\end{equation*}%
by the orthogonality relation \eqref{ortho-rel}.

We could have used similar methods as in \cite[Proposition 4.2]{abgrro17},
but we provide a more direct proof, based on the expression of the complex
Hermite polynomials in terms of Laguerre functions (\ref{ComplexHermite}).

\begin{proposition}
\label{aux-2-orth} Denote by $D_{R}$ the disc of radius $R$ centered at $0$.
It then holds 
\begin{equation}
\int_{D_{R}}H_{j,r}(z,\overline{z})\overline{H_{j\prime ,r\prime }(z,%
\overline{z})}e^{-\pi \left\vert z\right\vert ^{2}}dz=C_{j,r,j^{\prime
},r\prime }(R)\cdot {\delta }_{j-r-j\prime +r\prime }\text{,}
\label{doubleHermite}
\end{equation}%
with 
\begin{equation*}
C_{j,r,j^{\prime },r\prime }(R)=\sqrt{\frac{r!r\prime !}{j!j\prime !}}\pi
^{j-r}\int_{0}^{\pi R^{2}}L_{r}^{j-r}(t)L_{r\prime }^{j\prime -r\prime
}(t)e^{-t}dt.
\end{equation*}%
For $j=r$ and $j\prime =r\prime $ we obtain 
\begin{equation*}
\int_{D_{R}}H_{r,r}(z,\overline{z})\overline{H_{r\prime ,r\prime }(z,%
\overline{z})}e^{-\pi \left\vert z\right\vert ^{2}}dz=\int_{0}^{\pi R^{2}}{%
L_{r}^{0}(t)L_{r\prime }^{0}(t)}e^{-t}dt
\end{equation*}%
and, for $j=j\prime $,%
\begin{equation*}
\int_{D_{R}}H_{j,r}(z,\overline{z})\overline{H_{j,r\prime }(z,\overline{z})}%
e^{-\pi \left\vert z\right\vert ^{2}}dz=C_{j,r,r^{\prime }}(R)\cdot {\delta }%
_{r\prime -r}\text{.}
\end{equation*}
\end{proposition}

\textbf{Proof:\ }To avoid dividing the proof in two cases, we will use the
identity%
\begin{equation*}
\frac{(-x)^{j}}{j!}L_{r}^{j-r}(x)=\frac{(-x)^{r}}{r!}L_{j}^{r-j}(x)\text{,}
\end{equation*}%
to write (\ref{ComplexHermite}) as%
\begin{equation}
H_{j,r}(z,\overline{z})=\sqrt{\frac{r!}{j!}}\pi ^{\frac{j-r}{2}%
}z^{j-r}L_{r}^{j-r}\left( \pi \left\vert z\right\vert ^{2}\right) ,\qquad
j,r\in \mathbb{N}_{0}\text{.}  \label{ComplexHermite2}
\end{equation}%
Thus, 
\begin{align*}
& \int_{D_{R}}H_{j,r}(z,\overline{z})\overline{H_{j\prime ,r\prime }(z,%
\overline{z})}e^{-\pi \left\vert z\right\vert ^{2}}dz \\
& =\int_{D_{R}}{\sqrt{\frac{r!}{j!}}\pi ^{\frac{j-r}{2}}z^{j-r}L_{r}^{j-r}%
\left( \pi \left\vert z\right\vert ^{2}\right) }\overline{{\sqrt{\frac{%
r\prime !}{j\prime !}}\pi ^{\frac{j\prime -r\prime }{2}}z^{j\prime -r\prime
}L_{r\prime }^{j\prime -r\prime }\left( \pi \left\vert z\right\vert
^{2}\right) }}e^{-\pi \left\vert z\right\vert ^{2}}dz \\
& ={\sqrt{\frac{r!r\prime !}{j!j\prime !}}\pi ^{\frac{j-r+j\prime -r\prime }{%
2}}}\int_{D_{R}}{z^{j-r}}\overline{{z}}{^{j\prime -r\prime
}L_{r}^{j-r}\left( \pi \left\vert z\right\vert ^{2}\right) L_{r\prime
}^{j\prime -r\prime }\left( \pi \left\vert z\right\vert ^{2}\right) }e^{-\pi
\left\vert z\right\vert ^{2}}dz\text{.}
\end{align*}%
Setting $z=\rho e^{i\theta }$ we obtain%
\begin{align*}
& \int_{D_{R}}H_{j,r}(z,\overline{z})\overline{H_{j\prime ,r\prime }(z,%
\overline{z})}e^{-\pi \left\vert z\right\vert ^{2}}dz \\
& ={\sqrt{\frac{r!r\prime !}{j!j\prime !}}\pi ^{\frac{j-r+j\prime -r\prime }{%
2}}}\int_{0}^{R}\int_{0}^{2\pi }\rho ^{j-r-r^{\prime }+j\prime +1}e^{i\theta
(j-j\prime -r+r^{\prime })}{L_{r}^{j-r}\left( \pi \rho ^{2}\right)
L_{r\prime }^{j\prime -r\prime }\left( \pi \rho ^{2}\right) }e^{-\pi \rho
^{2}}d\rho d\theta \\
& ={\sqrt{\frac{r!r\prime !}{j!j\prime !}}\pi ^{\frac{j-r+j\prime -r\prime }{%
2}}}\int_{0}^{2\pi }e^{i\theta (j-j\prime -r+r^{\prime })}d\theta
\int_{0}^{R}\rho ^{j-r-r^{\prime }+j\prime +1}{L_{r}^{j-r}\left( \pi \rho
^{2}\right) L_{r\prime }^{j\prime -r\prime }\left( \pi \rho ^{2}\right) }%
e^{-\pi \rho ^{2}}d\rho \\
& ={\delta }_{j-j\prime -r+r^{\prime }}{\sqrt{\frac{r!r\prime !}{j!j\prime !}%
}\pi ^{j-r}\int_{0}^{R}2\pi \rho L_{r}^{j-r}(\pi \rho ^{2})L_{r\prime
}^{j\prime -r\prime }(\pi \rho ^{2})e^{-\pi \rho ^{2}}d\rho }\text{.}
\end{align*}%
\hfill $\Box $

Proposition \ref{aux-2-orth} now yields the local reproducing formula for
Hermite windows.

\begin{theorem}
\label{thm-double-orthog} For every $R>0$ and every $r,j\in \mathbb{N}_{0}$
one has 
\begin{equation}
V_{h_{r}}f(z)=C_{j,r}(R)^{-1}\int_{z+D_{R}}V_{h_{r}}f(w)K_{h_{j}}(z,w)dw%
\text{,}  \label{double-orthog}
\end{equation}%
with 
\begin{equation*}
C_{j,r}(R):=\langle \chi _{D_{R}}\cdot V_{h_{r}}{h_{r}},V_{h_{j}}h_{j}%
\rangle =\int_{0}^{\pi R^{2}}L_{r}^{0}(t)L_{j}^{0}(t)e^{-t}dt\text{.}
\end{equation*}
\end{theorem}

\textbf{Proof:\ }Rewriting Proposition \ref{aux-2-orth} using%
\begin{equation*}
V_{h_{r}}h_{j}(x,-\xi )=e^{i\pi x\xi -\tfrac{\pi }{2}\left\vert z\right\vert
^{2}}H_{j,r}(z,\bar{z})\text{,}
\end{equation*}%
leads to%
\begin{equation*}
\int_{D_{R}}V_{h_{r}}h_{j}(x,-\xi )\overline{V_{h_{r\prime }}h_{j\prime
}(x,-\xi )}dz=C_{j,r,j^{\prime },r\prime }(R)\cdot {\delta }_{j-r-j\prime
+r\prime }\text{.}
\end{equation*}%
In the case $j=r$\ we obtain by a change of variables 
\begin{equation*}
\int_{D_{R}}V_{h_{j}}h_{j}(z)\overline{V_{h_{j\prime }}h_{r\prime }(z)}%
dz=\int_{D_{R}}H_{j,j}(z,\overline{z})\overline{H_{j\prime ,r\prime }(z,%
\overline{z})}e^{-\pi \left\vert z\right\vert ^{2}}dz={\delta }_{j\prime
,r\prime }\cdot C_{j,j\prime }(R)\text{,}
\end{equation*}%
which implies that the following holds weakly 
\begin{equation*}
\int_{D_{R}}V_{h_{j}}h_{j}(w)\pi (w)h_{r}dw=C_{j,r}(R)\cdot h_{r}.
\end{equation*}%
Using \eqref{compos-tf} and \eqref{adjoint-tf} it thus follows 
\begin{align*}
V_{h_{r}}f(z)& =\langle f,\pi (z)h_{r}\rangle
=C_{j,r}(R)^{-1}\int_{D_{R}}\langle f,\pi (z)\pi (w)h_{r}\rangle \langle \pi
(w)h_{j},h_{j}\rangle dw \\
& =C_{j,r}(R)^{-1}\int_{D_{R}}e^{2\pi ix\eta }\langle f,\pi
(z+w)h_{r}\rangle \langle \pi (w)h_{j},h_{j}\rangle dw \\
& =C_{j,r}(R)^{-1}\int_{z+D_{R}}e^{2\pi ix(\eta -\xi )}\langle f,\pi
(w)h_{r}\rangle \langle \pi (w-z)h_{j},h_{j}\rangle dw \\
& =C_{j,r}(R)^{-1}\int_{z+D_{R}}\langle f,\pi (w)h_{r}\rangle \langle \pi
(w)h_{j},\pi (z)h_{j}\rangle dw\text{.}
\end{align*}%
\hfill $\Box $

Another consequence of Proposition \ref{aux-2-orth} is the following local
inversion formula, which allows to reconstruct $f$ from the values of the
STFT on arbitrary discs:

\begin{theorem}
\label{local-inversion} For every $R>0$, $r\in \mathbb{N}_{0}$ and $z\in 
\mathbb{R}^{2}$ we have 
\begin{equation}
f=\sum_{j\in \mathbb{N}_{0}}\Big(C_{j,r}(R)^{-1}\int_{z+D_{R}}V_{h_{r}}f(w)%
\langle \pi (w)h_{r},\pi (z)h_{j}\rangle dw\Big)\pi (z)h_{j}\text{,}\qquad
\forall f\in L^{2}(\mathbb{R})\text{.}  \label{eq-loc-inv}
\end{equation}
\end{theorem}

\textbf{Proof: }Write $f$ with respect to the orthonormal basis $%
\{h_{j}\}_{j\in \mathbb{N}_{0}}$ 
\begin{equation*}
f=\sum_{j\in \mathbb{N}_{0}}a_{j}h_{j}\text{.}
\end{equation*}%
By linearity of the STFT, one has%
\begin{equation*}
V_{h_{r}}f(w)=\sum_{j\in \mathbb{N}_{0}}a_{j}\langle h_{j},\pi
(w)h_{r}\rangle \text{.}
\end{equation*}%
Now, Proposition \ref{aux-2-orth} gives 
\begin{eqnarray*}
\int_{D_{R}}\langle f,\pi (w)h_{r}\rangle \langle \pi (w)h_{r},h_{k}\rangle
dw &=&\sum_{j\in \mathbb{N}_{0}}a_{j}\int_{D_{R}}\langle h_{j},\pi
(w)h_{r}\rangle \langle \pi (w)h_{r},h_{k}\rangle dw \\
&=&\sum_{j\in \mathbb{N}_{0}}a_{j}\int_{D_{R}}H_{j,r}(w,\overline{w})%
\overline{H_{k,r}(w,\overline{w})}e^{-\pi |w|^{2}}dw \\
&=&a_{k}C_{k,r}(R)\text{.}
\end{eqnarray*}%
Thus, 
\begin{equation*}
f=\sum_{j\in \mathbb{N}_{0}}a_{j}h_{j}=\sum_{j\in \mathbb{N}%
_{0}}C_{j,r}(R)^{-1}h_{j}\int_{D_{R}}\langle f,\pi (w)h_{r}\rangle \langle
\pi (w)h_{r},h_{j}\rangle dw.
\end{equation*}%
Applying this equality to $\pi (z)^{\ast }f$ and using the same argument as
in the proof of Theorem \ref{thm-double-orthog}, yields 
\begin{align*}
\pi (z)^{\ast }f& =\sum_{j\in \mathbb{N}_{0}}C_{j,r}(R)^{-1}h_{j}%
\int_{D_{R}}\langle \pi (z)^{\ast }f,\pi (w)h_{r}\rangle \langle \pi
(w)h_{r},h_{j}\rangle dw \\
& =\sum_{j\in \mathbb{N}_{0}}C_{j,r}(R)^{-1}h_{j}\int_{z+D_{R}}\langle f,\pi
(w)h_{r}\rangle \langle \pi (w)h_{r},\pi (z)h_{j}\rangle dw\text{.}
\end{align*}%
Now apply $\pi (z)$ on both sides to conclude the proof. \hfill $\Box $

\begin{remark}
\label{rem-form-Cn} If $j=r$, then $C_{r}(R)=1-e^{-\pi R^{2}}P_{r}(\pi
R^{2}) $, where $P_{r}$ is a polynomial of degree $2r$, $P_{r}(0)=1$ and $%
P_{0}\equiv 1$. See the appendix of \cite{Hutnik} for detailed calculations.
\end{remark}

\begin{remark}
\label{remark-poly} The so-called true (or pure, according to \cite%
{HendHaimi,HaiHen2}) polyanalytic Fock space $\mathcal{F}^{j}(\mathbb{C})$,
which can be defined as the span of $\{H_{j,r}(z,\overline{z})\}_{r\in 
\mathbb{N}}$ in $\mathcal{L}_2(\mathbb{C}):=L^2(\mathbb{C},e^{-\pi
\left\vert z\right\vert ^{2}})$ (see Sections 6.3 and 6.4 for more details)
or, equivalently, as the subspace of $\mathcal{L}_2(\mathbb{C})$ whose
elements satisfy the reproducing formula%
\begin{equation}
F(z)=\int_{\mathbb{C}}F(w)L_{j}^{0}(\pi \left\vert z-w\right\vert
^{2})e^{\pi z\overline{w}}e^{-\pi |w|^{2}}dw\text{,}  \label{RepKernel}
\end{equation}%
From Theorem \ref{thm-double-orthog} it follows that \eqref{doubleHermite}
implies the following local reproducing formula for $F\in \mathcal{F}^{j}(%
\mathbb{C})$: 
\begin{equation}
F(z)=C_{j,r}(R)^{-1}\int_{z+D_{R}}F(w)L_{r}^{0}(\pi \left\vert
z-w\right\vert ^{2})e^{\pi z\overline{w}}e^{-\pi |w|^{2}}dw\text{.}
\label{LocalPoly}
\end{equation}%
For $j=r$ this is what one would expect as a local version of (\ref%
{RepKernel}) and as an extension of the following local reproducing formula
for functions in the analytic Fock space $\mathcal{F}(\mathbb{C})=\mathcal{F}%
^{0}(\mathbb{C})$,\ obtained by Seip in \cite{Seip0} (this is also implicit
in \cite{Daubechies}): 
\begin{equation}
F(z)=(1-e^{-\pi R^{2}})^{-1}\int_{z+D_{R}}F(w)e^{\pi z\overline{w}}e^{-\pi
|w|^{2}}dw\text{.}  \label{local-reproducing}
\end{equation}%
However, for $j\neq r$ the spaces $\mathcal{F}^{j}(\mathbb{C})$ and $%
\mathcal{F}^{r}(\mathbb{C})$ are orthogonal and one could hardly expect %
\eqref{LocalPoly} to be true since, for every $F\in \mathcal{F}^{j}(\mathbb{C%
})$,%
\begin{equation*}
\int_{\mathbb{C}}F(w)L_{r}^{0}(\pi \left\vert z-w\right\vert ^{2})e^{\pi z%
\overline{w}}e^{-\pi |w|^{2}}dw=0\text{,}\qquad \forall z\in \mathbb{C}\text{%
.}
\end{equation*}%
Consequently, for $j\neq r$, the formula \eqref{LocalPoly} only holds for
finite $R$.
\end{remark}

\begin{remark}
It is not clear to us whether there exist other window functions that allow
for double orthogonality in sequences of concentric domains other than the
disc.
\end{remark}

\subsection{Estimates with explicit constants}

\label{sec:local}

We are now ready to formulate our main localization result for the Hermite
functions. Recall the maximum Nyquist density 
\begin{equation}
\rho (\Delta ,R):=\sup_{z\in \mathbb{R}^{2}}|\Delta \cap (z+D_{R})|.
\label{def-max-nyq}
\end{equation}%
We will also make use of the following notion of density 
\begin{equation}
A_{r}(\Delta ,R):=\sup_{z\in \mathbb{R}^{2}}\int_{\Delta \cap
z+D_{R}}|L_{r}^{0}(\pi |z-w|^{2})|e^{-\pi |z-w|^{2}/2}dw.  \label{Ar-density}
\end{equation}

\begin{theorem}
\label{main} Let $\Delta \subset \mathbb{R}^{2}$ and $f\in M^{p}$, $1\leq
p<\infty$. For every $0<R<\infty $, it holds 
\begin{equation}  \label{theo1-eq}
\dfrac{\Vert V_{h_r}f \cdot \chi_\Delta \Vert_{p}^p}{\Vert V_{h_r }f\Vert
_{p}^p}\leq \frac{A_{r} (\Delta,R)}{C_{r}(R)}\leq \frac{\rho(\Delta ,R)}{%
C_{r}(R)}\text{.} \ 
\end{equation}
\end{theorem}

\textbf{Proof:\ } In Proposition \ref{prop-measure-bound} take $%
K:=K_{h_{r}}\cdot \Omega _{R}$, where $\Omega _{R}(z,w):=\chi _{D_{R}}(z-w)$%
. Then $\theta (K)=1/C_{r}(R)$. Thus, if $d\mu (z)=\chi _{\Delta }dz$, we
have%
\begin{equation*}
\dfrac{\Vert V_{h_{r}}f\cdot \chi _{\Delta }\Vert _{1}}{\Vert
V_{h_{r}}f\Vert _{1}}\leq \frac{1}{C_{r}(R)}\Vert K_{h_{r}}\cdot \Omega
_{R}\Vert _{\mathcal{A}_{\chi _{\Delta }dz}}\text{.}
\end{equation*}%
Using the explicit formula \eqref{abs-kernel}, 
\begin{align*}
\Vert K_{h_{r}}\cdot \Omega _{R}\Vert _{\mathcal{A}_{\chi _{\Delta }dz}}&
=\sup_{z\in \mathbb{R}^{2}}\int_{\Delta }|K_{h_{r}}(z,w)|\chi _{D_{R}}(z-w)dw
\\
& =\sup_{z\in \mathbb{R}^{2}}\int_{\Delta \cap z+D_{R}}|L_{r}^{0}(\pi
|z-w|^{2})|e^{-\pi |z-w|^{2}/2}dw \\
& =A_{r}(\Delta ,R) \\
& \leq \rho (\Delta ,R).
\end{align*}%
Hence, the result holds for $p=1$. As $M^{1}\cap M^{\infty }=M^{1}$ is dense
in $M^{p}$ and 
\begin{equation*}
\sup_{\Vert f\Vert _{M^{\infty }}=1}\Vert V_{h_{r}}f\cdot \chi _{\Delta
}\Vert _{\infty }=1,
\end{equation*}%
the result for $1<p<\infty $ follows from Corollary \ref{measure-coroll}.
\hfill $\Box $

\begin{remark}
Results in the spirit of Theorem \ref{main} can be found for example in \cite%
[Section 4]{AnnihilatingSets} or \cite{Ascensi}. The estimates there are
however only given for sets with particular geometry, e.g. sets that are
thin at infinity or have finite Lebesgue measure, or without explicit
constants.
\end{remark}

An immediate consequence of Theorem \ref{main} is the following refined
(local) $L^{p}$-uncertainty principle for the short-time Fourier transform
(see \cite[Proposition 3.3.1]{Charly} and \cite{Jaming,groma13,rito14} for
other uncertainty principles for the STFT).

\begin{corollary}
\label{uncertainty} Suppose that $f\in M^{p}$, $1\leq p<\infty $, satisfies $%
\Vert V_{h_{r}}f\Vert _{p}=1$ and that $\Delta \subset \mathbb{R}^{2}$ and $%
\varepsilon \geq 0$ are such that 
\begin{equation*}
1-\varepsilon \leq \int_{\Delta }|V_{h_{r}}f(z)|^{p}dz\text{.}
\end{equation*}%
Then 
\begin{equation*}
1-\varepsilon \leq \inf_{R>0}\left( \frac{\rho (\Delta ,R)}{C_{r}(R)}\right)
\leq |\Delta |.
\end{equation*}
\end{corollary}

\begin{remark}
Essentially, Corollary \ref{uncertainty} states that the short-time Fourier
transform of a function in $M^{p}$ (using a Hermite window) cannot be well
concentrated on sets that are locally small over the entire time-frequency
plane. For an explicit example set $R=1$ and $r=0$. Then $C_{0}(R)=1-e^{-\pi
}\approx 0.96$, which implies that 
\begin{equation*}
\rho (\Delta ,1)\geq 0.95(1-\varepsilon ).
\end{equation*}%
Moreover, if we choose $\varepsilon =0.01$, there exists a subset of $\Delta 
$ contained in a disc of radius one, covering at least approximately $3/10$
of the area of that disc.
\end{remark}

\begin{remark}
If $K_{g}$ in (\ref{prop-general-g-eq}) shows sufficient off diagonal decay,
then the bound \eqref{prop-general-g-eq} behaves in a similar way as the
local integral $A_{r}(\Delta ,R)$ for the Hermite functions.
\end{remark}

\section{Optimality and sparse sets}

\label{sec:optimality}

Donoho and Logan \cite[Chapter 5]{DonohoLogan} discussed optimality of the
constant in \eqref{DLl1} as well as their $L^{2}$-estimate. As it turns out,
using extremal functions like the Beurling-Selberg function \cite{Beurling}
gives optimal constants within their method. In the STFT setup considered in
this paper, as far our knowledge goes, there is no theory of extremal
functions available. In the case of Gaussian window, we believe that our
local reproducing kernel is at least optimal among all kernels obtained from
truncating functions in $\mathcal{V}_{\varphi }$ on $D_{R}$ as $V_{\varphi
}\varphi $ optimizes the concentration problem on the disc for any $p\geq 1$.

Large sieve inequalities are particularly powerful if the localization
domain is sparse. It is nevertheless interesting to test the estimates in
cases where the solution of the localization problem is known. Note that
neither Donoho Logan's result applied to localization on an interval nor
Theorem \ref{main} applied to a disc achieve the actual solution. But this
is to be expected as the estimates hold for general sets. 

With a view to comparing the estimated and actual values
in cases where the exact solution is known, consider $\Delta =D_{R}$ and the
Gaussian window $g=\varphi =h_{0}$. In this case it is well known that the
Gaussian maximizes the concentration of the short-time Fourier transform in $%
D_{R}$ \cite{Daubechies, Seip0}. The $p$-norm can be explicitly evaluated as
follows: 
\begin{equation*}
\int_{D_{R}}|V_{\varphi }\varphi (z)|^{p}dz=\int_{D_{R}}e^{-\pi
p|z|^{2}/2}dz=2\pi \int_{0}^{R}\rho e^{-\pi p\rho ^{2}/2}d\rho =\frac{2}{p}%
(1-e^{-\pi pR^{2}/2})\text{.}
\end{equation*}%
Therefore, $\Vert V_{\varphi }\varphi \Vert _{p}^{p}=\frac{2}{p}$ and the
optimal solution of the concentration problem on $D_{\rho }$ is given by 
\begin{equation}
\sup_{f\in M^{p}}\frac{\int_{D_{R}}|V_{\varphi }f(z)|^{p}dz}{\Vert
V_{\varphi }f\Vert _{p}^{p}}=\frac{\int_{D_{R}}|V_{\varphi }\varphi
(z)|^{p}dz}{\Vert V_{\varphi }\varphi \Vert _{p}^{p}}=(1-e^{-\pi pR^{2}/2})%
\text{.}  \label{eq-gauss-disc}
\end{equation}%
Moreover, 
\begin{equation*}
A_{0}(D_{\rho },R)=\sup_{z\in \mathbb{R}^{2}}\int_{D_{\rho }\cap
z+D_{R}}e^{-\pi |z-w|^{2}/2}dw=\int_{D_{T}}e^{-\pi |w|^{2}/2}dw=2(1-e^{-\pi
T^{2}/2})\text{,}
\end{equation*}%
with $T=\min \{\rho ,R\}$ and $C_{0}(R)=(1-e^{-\pi R^{2}})$. Using Theorem %
\ref{main} we obtain the concentration estimate 
\begin{equation*}
\inf_{R>0}\frac{A_{0}(D_{\rho },R)}{C_{0}(R)}=\inf_{R>0}\frac{2(1-e^{-\pi
T^{2}/2})}{(1-e^{-\pi R^{2}})}=2(1-e^{-\pi \rho ^{2}/2})\text{.}
\end{equation*}%
Comparing our general estimate with the actual optimal value from %
\eqref{eq-gauss-disc} we observe that the estimate is not optimal for any $%
p\geq 1$. Let for example $p=1$. Then 
\begin{equation*}
\frac{2(1-e^{-\pi \rho ^{2}/2})}{1-e^{-\pi \rho ^{2}/2}}=2\text{.}
\end{equation*}%
For\ $p\in \lbrack 1,\infty \lbrack $ write $\mathcal{L}_{p}(%
%TCIMACRO{\U{2102} }%
%BeginExpansion
\mathbb{C}
%EndExpansion
)$ to denote the Banach space of all measurable functions equipped with the
norm 
\begin{equation*}
\left\Vert F\right\Vert _{\mathcal{L}_{p}(%
%TCIMACRO{\U{2102} }%
%BeginExpansion
\mathbb{C}
%EndExpansion
)}=\Big(\int_{%
%TCIMACRO{\U{2102} }%
%BeginExpansion
\mathbb{C}
%EndExpansion
}\left\vert F(z)\right\vert ^{p}e^{-\pi p\frac{\left\vert z\right\vert }{2}%
^{2}}\,dz\Big)^{1/p}\text{.}
\end{equation*}%
Now we turn our focus to the asymptotics of the concentration problem.
Define the distance of two sets in a standard way via 
\begin{equation*}
\,\mathrm{dist}(A,B):=\inf \{|x-y|:\ x\in A,\ y\in B\}\text{.}
\end{equation*}%
Let us consider the case where $\Delta $ is given by a finite union of sets $%
\Delta _{k}$ with increasing separation%
\begin{equation*}
d:=\min_{k\neq l}\,\mathrm{dist}(\Delta _{k},\Delta _{l})\text{.}
\end{equation*}%
It is easy to see that $A_{r}(\Delta ,R)\rightarrow \max_{k}A_{r}(\Delta
_{k},R)$ as $d\rightarrow \infty $. As we will show below, the concentration
problem is accurately described by this observation: it is decoupled. A
related result for the case of one dimensional band-limited functions was
derived in \cite[Theorem 10]{DonohoStark}.

\begin{proposition}
Let $g\in M^1$ and $\Delta$ be the union of $N$ disjoint, compact sets $%
\Delta_1,...,\Delta_N$. If $d$ tends to infinity, then 
\begin{equation}
\sup_{f\in L^2(\mathbb{R})}\frac{\|V_g f\cdot \chi_\Delta\|_2}{\|V_g f\|_2}%
\longrightarrow \max_{k=1,..,N}\sup_{f\in L^2(\mathbb{R})}\frac{\|V_g f\cdot
\chi_{\Delta_k}\|_2}{\|V_g f\|_2}.
\end{equation}
\end{proposition}

\textbf{Proof:\ } First, it trivially holds that 
\begin{equation*}
\max_{k=1,..,N}\sup_{f\in L^{2}(\mathbb{R})}\frac{\Vert V_{g}f\cdot \chi
_{\Delta _{k}}\Vert _{2}}{\Vert V_{g}f\Vert _{2}}\leq \sup_{f\in L^{2}(%
\mathbb{R})}\frac{\Vert V_{g}f\cdot \chi _{\Delta }\Vert _{2}}{\Vert
V_{g}f\Vert _{2}}.
\end{equation*}%
We can restrict ourselves to the case $\Delta =\Delta _{1}\cup \Delta _{2}$.
The general result then follows by induction. Now, for simplicity assume
that $\Vert g\Vert _{2}=1$ and define $f_{k},\ k=1,2,$ via 
\begin{equation*}
f_{k}:=V_{g}^{\ast }\big(V_{g}f\cdot \chi _{U_{d/3}(\Delta _{k})}\big),
\end{equation*}%
where $U_{d}(\Delta ):=\{z\in \mathbb{R}^{2}:\ \,\mathrm{dist}(z,\Delta
)\leq d\}$. Let $\varepsilon >0$ and choose $d=d(\varepsilon )$ large enough
such that for all $f\in L^{2}(\mathbb{R})$ and $k\in \{1,2\}$ 
\begin{equation}
\Vert \big(V_{g}f-V_{g}f_{k}\big)\cdot \chi _{\Delta _{k}}\Vert _{2}\leq
\varepsilon \Vert V_{g}f\Vert _{2},  \label{inside-bound}
\end{equation}%
and 
\begin{equation}
\Vert V_{g}f_{k}\cdot \chi _{U_{d/2}(\Delta _{k})^{c}}\Vert _{2}\leq
\varepsilon \Vert V_{g}f\Vert _{2}.  \label{outside-bound}
\end{equation}%
It is indeed possible to choose $d$ accordingly, since 
\begin{align*}
\Vert \big(V_{g}f-V_{g}f_{k}\big)\cdot \chi _{\Delta _{k}}\Vert _{2}^{2}&
=\int_{\Delta _{k}}\Big|V_{g}f(z)-\int_{U_{d/3}(\Delta
_{k})}V_{g}f(w)K_{g}(z,w)dw\Big|^{2}dz \\
& =\int_{\Delta _{k}}\Big|\int_{\mathbb{R}^{2}\backslash U_{d/3}(\Delta
_{k})}V_{g}f(w)K_{g}(z,w)dw\Big|^{2}dz \\
& \leq \int_{\Delta _{k}}\int_{\mathbb{R}^{2}\backslash U_{d/3}(\Delta
_{k})}|K_{g}(z,w)|^{2}dwdz\ \Vert V_{g}f\Vert _{2}^{2} \\
& \leq |\Delta _{k}|\sup_{z\in \Delta _{k}}\int_{\mathbb{R}^{2}\backslash
U_{d/3}(\Delta _{k})}|\langle g,\pi (z-w)g\rangle |^{2}dw\ \Vert V_{g}f\Vert
_{2}^{2} \\
& \leq |\Delta _{k}|\int_{\mathbb{R}^{2}\backslash D_{d/3}}|\langle g,\pi
(w)g\rangle |dw\ \Vert V_{g}f\Vert _{2}^{2} \\
& =C(k,d)\Vert V_{g}f\Vert _{2}^{2}\text{,}
\end{align*}%
where the last inequality follows from $|z-w|\geq d/3$, for $z\in \Delta
_{k} $ and $w\in \mathbb{R}^{2}\backslash U_{d/3}$. From the STFT being an
isometry it now follows that $C(k,d)\rightarrow 0$ as $d\rightarrow \infty $%
. To show that \eqref{outside-bound} is satisfied if $d$ is chosen big
enough, observe at first that 
\begin{equation*}
\sup_{z\in U_{d/2}(\Delta _{k})^{c}}|V_{g}f_{k}(z)|\leq \sup_{z\in
U_{d/2}(\Delta _{k})^{c}}\int_{U_{d/3}(\Delta
_{k})}|V_{g}f(w)K_{g}(z,w)|dw\leq \Vert V_{g}f\Vert _{\infty }\Vert
V_{g}g\Vert _{1}.
\end{equation*}%
The $L^{1}$-norm on the other hand can be estimated as 
\begin{align*}
\Vert V_{g}f_{k}\cdot \chi _{U_{d/2}(\Delta _{k})^{c}}\Vert _{1}& \leq
\int_{U_{d/2}(\Delta _{k})^{c}}\int_{U_{d/3}(\Delta
_{k})}|V_{g}f(w)K_{g}(z,w)|dwdz \\
& \leq \sup_{w\in U_{d/3}(\Delta _{k})}\int_{U_{d/2}(\Delta
_{k})^{c}}|K_{g}(z,w)|dz\ \Vert V_{g}f\Vert _{1} \\
& \leq \int_{\mathbb{R}^{2}\backslash D_{d/3}}|\langle g,\pi (z)g\rangle
|dz\ \Vert V_{g}f\Vert _{1} \\
& =\widetilde{C}(k,d)\Vert V_{g}f\Vert _{1}.
\end{align*}%
Hence, \eqref{outside-bound} follows by interpolation. As $\Vert
V_{g}f_{k}\Vert _{2}\leq \Vert V_{g}f\Vert _{2}$ we deduce from %
\eqref{inside-bound} that 
\begin{align*}
\Vert V_{g}f\cdot \chi _{\Delta }\Vert _{2}^{2} &=\Vert V_{g}f\cdot \chi
_{\Delta _{1}}\Vert _{2}^{2}+\Vert V_{g}f\cdot \chi _{\Delta _{2}}\Vert
_{2}^{2} \\
&\leq \Vert V_{g}f_{1}\cdot \chi _{\Delta _{1}}\Vert _{2}^{2}+\Vert
V_{g}f_{2}\cdot \chi _{\Delta _{2}}\Vert _{2}^{2}+C\varepsilon \Vert
V_{g}f\Vert _{2}^{2}\text{.}
\end{align*}%
Moreover, by \eqref{outside-bound}, we have the following almost
orthogonality relation for $V_{g}f_{1}$ and $V_{g}f_{2}$: 
\begin{align*}
\Vert V_{g}f_{1}+V_{g}f_{2}\Vert _{2}^{2} &\geq \Vert V_{g}f_{1}\Vert
_{2}^{2}+\Vert V_{g}f_{2}\Vert _{2}^{2}-2|\langle
V_{g}f_{1},V_{g}f_{2}\rangle | \\
&\geq \Vert V_{g}f_{1}\Vert _{2}^{2}+\Vert V_{g}f_{2}\Vert
_{2}^{2}-2|\langle V_{g}f_{1}\cdot \chi _{U_{d/2}(\Delta
_{1})^{c}},V_{g}f_{2}\rangle |-2|\langle V_{g}f_{1},V_{g}f_{2}\cdot \chi
_{U_{d/2}(\Delta _{2})^{c}}\rangle | \\
&\geq \Vert V_{g}f_{1}\Vert _{2}^{2}+\Vert V_{g}f_{2}\Vert
_{2}^{2}-4\varepsilon \Vert V_{g}f_{1}\Vert _{2}\Vert V_{g}f_{2}\Vert _{2} \\
&\geq \Vert V_{g}f_{1}\Vert _{2}^{2}+\Vert V_{g}f_{2}\Vert
_{2}^{2}-4\varepsilon \Vert V_{g}f\Vert _{2}^{2}.
\end{align*}%
Now, as $\Vert V_{g}f_{1}+V_{g}f_{2}\Vert _{2}^{2}\leq \Vert V_{g}f\cdot
\chi _{\Delta _{1}\cup \Delta _{2}}\Vert _{2}^{2}\leq \Vert V_{g}f\Vert
_{2}^{2}$ it follows that 
\begin{align*}
\frac{\Vert V_{g}f\cdot \chi _{\Delta }\Vert _{2}^{2}}{\Vert V_{g}f\Vert
_{2}^{2}}& \leq (1+C\varepsilon )\frac{\Vert V_{g}f_{1}\cdot \chi _{\Delta
_{1}}\Vert _{2}^{2}+\Vert V_{g}f_{2}\cdot \chi _{\Delta _{2}}\Vert _{2}^{2}}{%
\Vert V_{g}f_{1}\Vert _{2}^{2}+\Vert V_{g}f_{2}\Vert _{2}^{2}}+C\varepsilon
\\
& \leq (1+C\varepsilon )\max_{k=1,2}\frac{\Vert V_{g}f_{k}\cdot \chi
_{\Delta _{k}}\Vert _{2}^{2}}{\Vert V_{g}f_{k}\Vert _{2}^{2}}+C\varepsilon \\
& \leq (1+C\varepsilon )\max_{k=1,2}\sup_{f\in L^{2}(\mathbb{R})}\frac{\Vert
V_{g}f\cdot \chi _{\Delta _{k}}\Vert _{2}^{2}}{\Vert V_{g}f\Vert _{2}^{2}}%
+C\varepsilon ,
\end{align*}%
which concludes the proof if we take the supremum over $L^{2}(\mathbb{R})$
on the left hand side. \hfill $\Box $

\begin{remark}
Although we expect a similar result to hold for the concentration problem in 
$M^p$ we were not able to prove it. The main problem is that our argument
relies on $\|V_gV_g^\ast\|_{2\rightarrow 2}=1$ which is not true on $M^p$, $%
p\neq 2$.
\end{remark}

Finally, we present a conjecture on an extremal problem of localization of
the STFT with Gaussian window which is the joint time-frequency analogue of 
\cite[Conjecture 1]{DonohoStark}.

\begin{conjecture}
\label{extremal-conjecture} Let $\Delta \subset \mathbb{R}^{2}$ be a set of
finite measure and $\varphi =h_{0}$ be the Gaussian. Then 
\begin{equation*}
\sup_{|\Delta |=A}\sup_{f\in M^{p}}\frac{\Vert V_{\varphi }f\cdot \chi
_{\Delta }\Vert _{p}^{p}}{\Vert V_{\varphi }f\Vert _{p}^{p}}
\end{equation*}%
is attained if and only if $\Delta =z+D_{\sqrt{A/\pi }}$ for some $z\in 
\mathbb{R}^{2}$, up to perturbations of Lebesgue measure zero.
\end{conjecture}

The next Proposition provides extra support for the conjecture, by showing
that the disc is the unique solution (up to perturbations of Lebesgue
measure zero) of a certain extremal problem. This will in turn imply that
the disc maximizes $A_{0}(\Delta ,R)$ for all $R>0$, where the area of $%
\Delta$ is fixed. Consequently, Conjecture 1 is backed by the estimates of
Theorem \ref{main}.

\begin{proposition}
Let $\alpha >0$. The disc $D_{R}$, $R=\sqrt{A/\pi }$ is the unique (up to
perturbations of Lebesgue measure zero) minimizer of the following extremal
problem: 
\begin{equation}
\sup_{\Omega \subset \mathbb{R}^{n}}\int_{\Omega }e^{-\alpha
|z|^{2}}dz,\qquad \mbox{subject to }|\Omega |=A.  \label{extrem-probl}
\end{equation}
\end{proposition}

\textbf{Proof:\ } Let us assume to the contrary that there exists $\Omega
\subset \mathbb{R}^{n}$, such that $|\Omega |=C$ and $|\Omega \backslash
D_{R}|>0$ which maximizes \eqref{extrem-probl}. Define $\Omega _{r}:=\Omega
\backslash D_{r}$. Then there exists $\varepsilon =\varepsilon (R,|\Omega
\backslash D_{R}|)>0$ such that $|\Omega _{R+\varepsilon }|\geq \frac{%
|\Omega _{R}|}{2}>0$. Let $I\subset D_{R}\backslash \Omega $ be any set that
satisfies $|I|=|\Omega _{R+\varepsilon }|$. (Such a set exists as $\Omega $
contains a set of size $|\Omega _{R+\varepsilon }|$ outside the disc $D_{R}$
and has the same size as the disc.) Define another set $\Omega ^{\ast
}:=\Omega \backslash \Omega _{R+\varepsilon }\cup I$. It then holds that $%
|\Omega ^{\ast }|=C$ and 
\begin{align*}
\int_{\Omega ^{\ast }}e^{-\alpha |z|^{2}}dz &=\int_{\Omega \backslash
D_{R+\varepsilon }}\hspace{-0.2cm}e^{-\alpha |z|^{2}}dz+\int_{I}e^{-\alpha
|z|^{2}}dz \\
&\geq \int_{\Omega \backslash \Omega _{R+\varepsilon }}\hspace{-0.2cm}%
e^{-\alpha |z|^{2}}dz+e^{-\alpha R^{2}}|I| \\
&>\int_{\Omega \backslash \Omega _{R+\varepsilon }}\hspace{-0.2cm}e^{-\alpha
|z|^{2}}dz+e^{-\alpha (R+\varepsilon )^{2}}|\Omega _{R+\varepsilon }| \\
&\geq \int_{\Omega \backslash \Omega _{R+\varepsilon }}\hspace{-0.2cm}%
e^{-\alpha |z|^{2}}dz+\int_{\Omega _{R+\varepsilon }}e^{-\alpha |z|^{2}}dz \\
&=\int_{\Omega }e^{-\alpha |z|^{2}}dz,
\end{align*}%
which contradicts the assumption that $\Omega $ maximizes %
\eqref{extrem-probl}. \hfill $\Box $

\subsection{Recovery of STFT measurements}

Now we will rephrase Proposition \ref{l1-sparse} and \ref{l1-missing} in the
context of reconstructing STFT-data using Theorem \ref{main}. Let $B_{1}=%
\mathcal{V}_{h_{r}}^{1}=V_{h_{r}}(M^{1})\subset L^{1}(\mathbb{R}^{2})$. By
the correspondence principle \eqref{corr-princ} we can replace minimization
on $B_{1}$ by minimization on $M^{1}$ (which is independent of the
particular choice of the order of the Hermite window).

\begin{corollary}
\label{sparsenoisecor} Suppose that $G=V_{h_r }f+N$ is observed, where $f\in
M^{1}$, $N\in L^{1}(\mathbb{R}^{2})$ and that the unknown support $\Delta $
of $N$ satisfies 
\begin{equation}
A_r(\Delta ,R)<\frac{C_{r}(R)}{2}\text{,}  \label{planarsparse}
\end{equation}%
for some $R>0$. Then $\delta (\Delta )<\frac{1}{2}$ and the solution of the
minimization problem 
\begin{equation*}
\beta (G)=\arg \min_{g\in M^{1}}\big\Vert G-V_{h_r }g\big\Vert_{1}
\end{equation*}%
is unique and recovers the signal $f$ perfectly ($\beta (G)=f$).
\end{corollary}

\begin{corollary}
\label{missingcorollary} Let $f\in M^{1}$ and suppose that one observes $%
H=P_{\Delta ^{c}}(V_{h_{r}}f+N)$, where $\Vert N\Vert _{1}\leq \varepsilon $
and that the domain $\Delta $ of missing data satisfies 
\begin{equation}
A_{r}(\Delta ,R)<C_{r}(R)\text{,}  \label{planarsparselessone}
\end{equation}%
for some $R>0$. Then any solution of 
\begin{equation*}
\sigma (H)=\arg \min_{g\in M^{1}}\big\Vert P_{\Delta ^{c}}(H-V_{h_{r}}g)%
\big\Vert_{1}
\end{equation*}%
satisfies 
\begin{equation*}
\left\Vert V_{h_{r}}\big(f-\sigma (H)\big)\right\Vert _{1}\leq \frac{%
2\varepsilon \cdot C_{r}(R)}{C_{r}(R)-A_{r}(\Delta ,R)}.
\end{equation*}
\end{corollary}

For other approaches to the recovery of sparse time-frequency
representations which concentrate on the set-up of finite sparse
time-frequency representations, see \cite{SparsityTF,PRT}.

\section{Extensions to other settings}

\label{sec:other-settings}

\subsection{Discrete Gabor systems}

In this section, we will apply Corollary \ref{measure-coroll} to discrete
Gabor systems. Let $\Lambda \subset \mathbb{R}^{2}$ be discrete and $\Delta
\subset \Lambda $. We define the discrete maximum Nyquist density by 
\begin{equation*}
\rho ^{d}(\Delta ,R):=\sup_{z\in \mathbb{R}^{2}}\#\{\Delta \cap z+D_{R}\}
\end{equation*}%
and consider also 
\begin{equation*}
A_{r}^{d}(\Delta ,R):=\sup_{z\in \mathbb{R}^{2}}\sum_{\lambda \in \Delta
\cap z+D_{R}}L_{r}^{0}(\pi |z-\lambda |^{2})e^{-\pi |z-\lambda |^{2}/2}.
\end{equation*}%
We define the measure $\mu $ to be $d\mu (\lambda ):=\delta _{\Lambda
}(\lambda )\cdot \chi _{\Delta }(\lambda )d\lambda $. First, observe that 
\begin{equation*}
\sup_{f\in M^{\infty }}\frac{\Vert V_{h_{r}}f|_{\Delta }\Vert _{\infty }}{%
\Vert V_{h_{r}}f\Vert _{\infty }}=1.
\end{equation*}%
Applying Corollary \ref{measure-coroll} yields 
\begin{equation*}
\frac{\sum_{\lambda \in \Delta }|V_{h_{r}}f(\lambda )|^{p}}{\Vert
V_{h_{r}}f\Vert _{p}^{p}}\leq \frac{A_{r}^{d}(\Delta ,R)}{C_{r}(R)}\leq 
\frac{\rho ^{d}(\Delta ,R)}{C_{r}(R)}\text{.}
\end{equation*}%
If we also assume that $\{\pi (\lambda )h_{r}\}_{\lambda \in \Lambda }$ is a
frame for $L^{2}(\mathbb{R})$, then, by Remark \ref{banach-gabor-frame},
there exists a lower Banach frame bound $m=m(r,\Lambda ,p)$ such that 
\begin{equation*}
m\Vert f\Vert _{M^{p}}^{p}=m\Vert V_{h_{r}}f\Vert _{p}^{p}\leq \sum_{\lambda
\in \Lambda }|V_{h_{r}}f(\lambda )|^{p}.
\end{equation*}%
It then holds 
\begin{equation}
\frac{\sum_{\lambda \in \Delta }|V_{h_{r}}f(\lambda )|^{p}}{\sum_{\lambda
\in \Lambda }|V_{h_{r}}f(\lambda )|^{p}}\leq \frac{A_{r}^{d}(\Delta ,R)}{%
m\cdot C_{r}(R)}\leq \frac{\rho ^{d}(\Delta ,R)}{m\cdot C_{r}(R)}\text{.}
\end{equation}%
If we would like to use this estimate for discrete signal recovery, then the
bound $\frac{A_{r}^{d}(\Delta ,R)}{m\cdot C_{r}(R)}$ should be small, or at
least less than one half. If $\Lambda $ is a lattice in $\mathbb{R}^{2}$
which does not deviate too much from the square lattice then $m$ scales with
the density of $\Lambda $. Since also $\rho ^{d}(\Delta ,R)$ and $%
A_{r}^{d}(\Delta ,R)$ show similar behavior, it is still possible to get
small concentration bounds if the density of $\Lambda $ is increased. An
interesting direction for further research could also be to study the
concentration problem for complete Gabor systems that are not frames (this
is the case of several lattice configurations for Gabor systems with Hermite
functions, which are known to be complete \cite{Gabor} but not frames \cite%
{Jacob}), but still allow for reconstruction using dual systems \cite%
{SpeckbacherBalazs}.

\subsection{Vector-valued STFT transforms}

Vector-valued time-frequency analysis \cite{DL} is motivated by the problem
of multiplexing of signals, where one wants to transmit several signals over
a single channel followed by separating and recovering the signals at the
receiver \cite{Balan}. A classical way to do this is to store the
information for every function in mutually orthogonal subspaces. The
orthogonality relation \eqref{ortho-rel} for the short-time Fourier
transform suggested a vector-valued version of the STFT using mutually
orthogonal windows, called the super Gabor transform \cite{Abreustructure},
in a reference to the connection with \cite{Fuhr, CharlyYurasuper}. In the
case of a vector constituted by Hermite functions, this reads 
\begin{equation*}
\mathbf{V}_{\mathbf{h_{n}}}\mathbf{f}%
(z)=V_{(h_{0},...,h_{n})}(f_{0},...,f_{n})(z):=%
\sum_{k=0}^{n}V_{h_{k}}f_{k}(z)\text{,}
\end{equation*}%
where $\mathbf{f}:=(f_{0},f_{1},...,f_{n})\in L^{2}(\mathbb{R})^{n+1}$ and $%
\mathbf{h_{n}}:=(h_{0},h_{1},...,h_{n})$ is the vector of the first $n+1$
Hermite functions. This is the continuous transform associated to Gabor
superframes with Hermite windows \cite{Fuhr, CharlyYurasuper} and to
sampling in polyanalytic Fock spaces \cite{Abreusampling}. 
% for the purpose of processing simultaneously $n$ signals using a vectorial window constituted by the first  $n$ Hermite functions. 
The function $f_{k}$ can then be reconstructed by 
\begin{equation*}
f_{k}=V_{h_{k}}^{\ast }\mathbf{V_{h_{n}}f}.
\end{equation*}%
The range of this transform is a Hilbert space with reproducing kernel given
by 
\begin{equation*}
\mathbf{K_{h_{n}}}(z,w)=e^{i\pi (x+y)(\omega -\eta )}L_{n}^{1}(\pi
|z-w|^{2})e^{-\pi |z-w|^{2}/2},
\end{equation*}%
(this follows from (\ref{RepKernel}) and the summation relation $%
\sum_{k=0}^{n}L_{k}^{\alpha }=L_{n}^{\alpha +1}$ of the Laguerre functions). 
%Consequently, we have  \begin{equation*}|\bm{K_{h_n}}(z,w)|=L_{n}^{1}(\pi |z-w|^{2})e^{-\pi |z-w|^{2}/2}.\end{equation*}
For the basis functions of the reproducing kernel $\mathbf{K_{h_{n}}}$,
double orthogonality is lost, since the cross terms are not zero (the
Laguerre functions are not orthogonal on any interval $[0,R]$, for $R<\infty 
$). We can however still define a local kernel that yields an estimate in
terms of $A_{0}(\Delta ,R)$.

%Of course we can also derive a similar result to Proposition \ref{loc-general-g} for $\mathbf{K_{h_n}}$ and obtain a concentration estimate since $\mathbf{K_{h_n}}\in L^1(\R^2)$.

Set $\mathcal{V}_{g}^{p}:=\{V_{g}f:\ f\in M^{p}\}$ then the orthogonal
decomposition extends to the modulation spaces $M^{p}$, see \cite{abgroe12}.
For $1\leq p<\infty $, we have 
\begin{equation*}
\mathbf{V_{h_{n}}}\Big(\prod_{k=0}^{n}M^{p}\Big)=\mathcal{V}%
_{h_{0}}^{p}\oplus \mathcal{V}_{h_{1}}^{p}\oplus ...\oplus \mathcal{V}%
_{h_{n}}^{p}.
\end{equation*}%
Therefore, 
\begin{equation*}
\Vert \mathbf{V_{h_{n}}f}\Vert _{p}=\Big\|\sum_{k=0}^{n}V_{h_{k}}f_{k}\Big\|%
_{p}\asymp \sum_{k=0}^{n}\Vert V_{h_{k}}f_{k}\Vert _{p}.
\end{equation*}%
It follows from Theorem \ref{thm-double-orthog} that 
\begin{equation*}
\mathbf{V_{h_{n}}f}(z)=\sum_{k=0}^{n}V_{h_{k}}f_{k}(z)=%
\sum_{k=0}^{n}C_{k,0}(R)^{-1}\int_{z+D_{R}}V_{h_{k}}f_{k}(w)K_{h_{0}}(z,w)dw,
\end{equation*}%
which yields 
\begin{align*}
\Vert \mathbf{V_{h_{n}}f}\cdot \chi _{\Delta }\Vert _{1}& \leq \max_{0\leq
m\leq n}C_{m,0}(R)^{-1}\int_{\Delta
}\sum_{k=0}^{n}\int_{z+D_{R}}|V_{h_{k}}f_{k}(w)K_{h_{0}}(z,w)|dwdz \\
& \leq \max_{0\leq m\leq n}C_{m,0}(R)^{-1}\cdot A_{0}(\Delta ,R)\cdot
\sum_{k=0}^{n}\Vert V_{h_{k}}f_{k}\Vert _{1} \\
& \leq \widetilde{C}\cdot \max_{0\leq m\leq n}C_{m,0}(R)^{-1}\cdot
A_{0}(\Delta ,R)\cdot \Big\|\sum_{k=0}^{n}V_{h_{k}}f_{k}\Big\|_{1} \\
& =\widetilde{C}\cdot \max_{0\leq m\leq n}C_{m,0}(R)^{-1}\cdot A_{0}(\Delta
,R)\cdot \Vert \mathbf{V_{h_{n}}f}\Vert _{1}\text{.}
\end{align*}%
Hence, we have shown that the concentration operator of a multiplexed
short-time Fourier transform can also be estimated in terms of $A_{0}(\Delta
,R)$ and $\rho (\Delta ,R)$ at the cost of a larger normalization constant
and the additional factor $\widetilde{C}$.

\subsection{True polyanalytic Fock spaces}

The Bargmann transform $\mathcal{B}$, defined as 
\begin{equation*}
\mathcal{B}f(z)=2^{\frac{1}{4}}\int_{%
%TCIMACRO{\U{211d} }%
%BeginExpansion
\mathbb{R}
%EndExpansion
}f(t)e^{2\pi tz-\pi z^{2}-\frac{\pi }{2}t^{2}}dt\text{,}
\end{equation*}%
is an isomorphism $\mathcal{B}:L^{2}(%
%TCIMACRO{\U{211d} }%
%BeginExpansion
\mathbb{R}
%EndExpansion
)\rightarrow \mathcal{F}_{2}(%
%TCIMACRO{\U{2102} }%
%BeginExpansion
\mathbb{C}
%EndExpansion
)$, where $\mathcal{F}_{2}(%
%TCIMACRO{\U{2102} }%
%BeginExpansion
\mathbb{C}
%EndExpansion
)$ is the classical Bargmann-Fock space of entire functions. One can define
a sequence of transforms $\mathcal{B}^{r+1}:L^{2}(%
%TCIMACRO{\U{211d} }%
%BeginExpansion
\mathbb{R}
%EndExpansion
)\rightarrow \mathcal{F}_{2}^{r+1}(%
%TCIMACRO{\U{2102} }%
%BeginExpansion
\mathbb{C}
%EndExpansion
)$ as a Hilbert space isomorphism mapping onto true polyanalytic Fock spaces 
\cite{Abreusampling,VasiFock} as follows: 
\begin{equation*}
\mathcal{B}^{r+1}f(z)=\left( \frac{\pi ^{r}}{r!}\right) ^{\frac{1}{2}}e^{\pi
\left\vert z\right\vert ^{2}}\left( \partial _{z}\right) ^{r}\left[ e^{-\pi
\left\vert z\right\vert ^{2}}\mathcal{B}f(z)\right]
\end{equation*}%
The relation between Gabor transforms with Hermite functions and true
polyanalytic Barg- mann transforms of general order $r$ reads \cite%
{Abreusampling}: 
\begin{equation}
e^{-i\pi x\xi +\pi \frac{\left\vert z\right\vert ^{2}}{2}}V_{h_{r}}f(x,-\xi
)=\mathcal{B}^{r+1}f(z)\text{.}  \label{Barg_hermite}
\end{equation}%
The $L^{p}$ version of the polyanalytic Bargmann-Fock spaces has been
introduced in \cite{abgroe12}, where the link to Gabor analysis has been
particularly useful. For\ $p\in \lbrack 1,\infty \lbrack $ write $\mathcal{L}%
_{p}(%
%TCIMACRO{\U{2102} }%
%BeginExpansion
\mathbb{C}
%EndExpansion
)$ to denote the Banach space of all measurable functions equipped with the
norm 
\begin{equation*}
\left\Vert F\right\Vert _{\mathcal{L}_{p}(%
%TCIMACRO{\U{2102} }%
%BeginExpansion
\mathbb{C}
%EndExpansion
)}=\Big(\int_{%
%TCIMACRO{\U{2102} }%
%BeginExpansion
\mathbb{C}
%EndExpansion
}\left\vert F(z)\right\vert ^{p}e^{-\pi p\frac{\left\vert z\right\vert }{2}%
^{2}}\,dz\Big)^{1/p}\text{.}
\end{equation*}%
As a corollary of Theorem \label{main copy(1)}, we thus obtain the inequality%
\begin{equation*}
\dfrac{\Vert F\cdot \chi _{\Delta }\Vert _{\mathcal{L}_{p}(%
%TCIMACRO{\U{2102} }%
%BeginExpansion
\mathbb{C}
%EndExpansion
)}^{p}}{\Vert F\Vert _{\mathcal{L}_{p}(%
%TCIMACRO{\U{2102} }%
%BeginExpansion
\mathbb{C}
%EndExpansion
)}^{p}}\leq \frac{\rho (\Delta ,R)}{C_{r}(R)},\qquad \forall F\in\mathcal{F}%
^{r+1}_p\text{.}
\end{equation*}

\subsection{Polyanalytic Fock spaces}

A function $F(z,\overline{z}),$ defined on a subset of $\mathbb{C}$, and
satisfying the generalized Cauchy-Riemann equations 
\begin{equation}
\left( \partial _{\overline{z}}\right) ^{n}F(z,\overline{z})=\frac{1}{2^{n}}%
\left( \frac{\partial }{\partial x}+i\frac{\partial }{\partial \xi }\right)
^{n}F(x+i\xi ,x-i\xi )=0\text{,}  \label{eq:c1}
\end{equation}%
is said to be \emph{polyanalytic of order }$n-1$.

\begin{definition}
We say that a function $F$ belongs to the \emph{polyanalytic Fock space} $%
\mathbf{F}_{2}^{n+1}\left( 
%TCIMACRO{\U{2102} }%
%BeginExpansion
\mathbb{C}
%EndExpansion
\right) $, if $\left\Vert F\right\Vert _{\mathcal{L}_{2}(\mathbb{C})}<\infty 
$ and $F$ is polyanalytic of order $n$.
\end{definition}

Polyanalytic Fock spaces seem to have been first considered by Balk \cite[%
pag. 170]{Balk}. Vasilevski \cite{VasiFock} obtained the following
decompositions in terms of the spaces $\mathcal{F}_{2}^{r}(%
%TCIMACRO{\U{2102} }%
%BeginExpansion
\mathbb{C}
%EndExpansion
)$:%
\begin{equation}
\mathbf{F}_{2}^{n}(%
%TCIMACRO{\U{2102} }%
%BeginExpansion
\mathbb{C}
%EndExpansion
)=\mathcal{F}_{2}^{1}(%
%TCIMACRO{\U{2102} }%
%BeginExpansion
\mathbb{C}
%EndExpansion
)\oplus ...\oplus \mathcal{F}_{2}^{n}(%
%TCIMACRO{\U{2102} }%
%BeginExpansion
\mathbb{C}
%EndExpansion
)  \label{orthogonal}
\end{equation}%
and%
\begin{equation*}
\mathcal{L}_{2}(\mathbb{C})=\bigoplus_{n=1}^{\infty }\mathcal{F}_{2}^{n}(%
%TCIMACRO{\U{2102} }%
%BeginExpansion
\mathbb{C}
%EndExpansion
).
\end{equation*}%
We can rewrite the transform of the previous section as a transform $\mathbf{%
B}^{n}:L^{2}(%
%TCIMACRO{\U{211d} }%
%BeginExpansion
\mathbb{R}
%EndExpansion
,%
%TCIMACRO{\U{2102} }%
%BeginExpansion
\mathbb{C}
%EndExpansion
^{n})\rightarrow \mathbf{F}^{n}(%
%TCIMACRO{\U{2102} }%
%BeginExpansion
\mathbb{C}
%EndExpansion
)$ mapping each vector $\mathbf{f=(}f_{1},...,f_{n}\mathbf{)}\in L^{2}(%
%TCIMACRO{\U{211d} }%
%BeginExpansion
\mathbb{R}
%EndExpansion
,%
%TCIMACRO{\U{2102} }%
%BeginExpansion
\mathbb{C}
%EndExpansion
^{n})$ to 
\begin{equation}
\mathbf{B}^{n}\mathbf{f=}e^{-i\pi x\xi +\pi \frac{\left\vert z\right\vert
^{2}}{2}}\mathbf{V}_{\mathbf{h}_{n-1}}\mathbf{f}(\lambda ).
\label{polyvectorBarg}
\end{equation}%
Since the multiplier $e^{-i\pi x\xi +\pi \frac{\left\vert z\right\vert ^{2}}{%
2}}$ in (\ref{Barg_hermite}) is the same for every $n$, we have: 
\begin{equation}
\mathbf{B}^{n}\mathbf{f=}\mathcal{B}^{1}f_{1}+...+\mathcal{B}^{n}f_{n}\,.
\label{polyanalyticBargmann}
\end{equation}%
This map is again a Hilbert space isomorphism and is called the \emph{%
polyanalytic Bargmann transform }\cite{Abreusampling}. The identity 
\begin{equation*}
\mathbf{V_{h_{n}}f}(z)=\sum_{k=0}^{n}C_{k,0}(R)^{-1}%
\int_{z+D_{R}}V_{h_{k}}f_{k}(w)K_{h_{0}}(z,w)dw
\end{equation*}%
can be written as%
\begin{equation*}
\mathbf{B}^{n}\mathbf{f}=\sum_{k=0}^{n}C_{k,0}(R)^{-1}\int_{z+D_{R}}\mathcal{%
B}^{k}f_{k}(w)K_{h_{0}}(z,w)e^{-\pi \left\vert w\right\vert ^{2}}dw\text{.}
\end{equation*}%
Rephrasing the discussion in the end of section 6.2, leads to the inequality 
\begin{equation*}
\Vert \mathbf{B}^{n}\mathbf{f}\cdot \chi _{\Delta }\Vert _{\mathcal{L}_{1}(%
%TCIMACRO{\U{2102} }%
%BeginExpansion
\mathbb{C}
%EndExpansion
)}\leq \widetilde{C}\cdot \max_{0\leq m\leq n}C_{m,0}(R)^{-1}\cdot
A_{0}(\Delta ,R)\cdot \Big\|\mathbf{B}^{n}\mathbf{f}\Big\|_{\mathcal{L}_1(%
\mathbb{C})}\text{.}
\end{equation*}

\section{Further questions}

\begin{enumerate}
\item If there exists a function $g\in L^{2}({\mathbb{R}})$ that allows for
a local reproducing formula on all discs of radius $R>0$, i.e., if 
\begin{equation*}
V_{g}f(z)=C_g(R)^{-1}\int_{z+D_{R}}\langle f,\pi (w)g\rangle \langle \pi
(w)g,\pi (z)g\rangle dw,\qquad \forall f\in L^{2}(\mathbb{R}),
\end{equation*}%
does it follow that $g$ is necessarily a Hermite function?

\item This problem concerns a generalization of the main result in \cite%
{abdoe12} using Hermite window instead of Gaussian window. If $\Omega $ is
simply connected and $h_{j}$ is an eigenfunction of the following
localization operator 
\begin{equation*}
H_{\Omega }^{r}f:=\int_{\Omega }\langle f,\pi (w)h_{r}\rangle \pi (w)h_{r}dw,
\end{equation*}%
does it follow that $\Omega $ is a disc centered at the origin?

\item Is it possible to find a window $g$ such that double orthogonality
holds in a sequence of non-circular domains $\Omega _{1}\subset \Omega
_{2}\subset ...\Omega _{\infty }=\mathbb{R}^{2}$?

\item Due to the orthogonality in concentric domains, the analysis in the
case of Hermite windows avoided the use of the extremal functions required,
for instance in \cite{DonohoLogan}. However, if one aims to extend the
results of \cite{DonohoLogan} to the challenging setup of general de Branges
spaces, used in the characterization of Fourier frames in \cite%
{FourierFrames}, such a simplification is unlikely to occur. It is thus a
natural question to ask if the results in \cite{ExtremaldeBranges} can be
used for this purpose. A related setup where one can expect the aid of
explicit formulas is the one of the band-limited multidimensional Fourier
transform of radial functions \cite{SlepianIV}, which essentially boils down
to the band-limited Hankel transform, where the localization operators and
the Nyquist rate have been studied in detail \cite{AB}.

\item Prove or disprove Conjecture \ref{extremal-conjecture}.
\end{enumerate}

\section*{Acknowledgement}

This work was funded by the Austrian Science Fund (FWF) START-project FLAME
('Frames and Linear Operators for Acoustical Modeling and Parameter
Estimation'; Y 551-N13) and FWF project `Operators and Time-Frequency
Analysis' P 31225-N32.

\bibliographystyle{plain}
\bibliography{paperbib}

\end{document}